\documentclass[11pt, a4paper]{article}
\usepackage{amsmath} \usepackage{euscript}
\usepackage{mymatrix2}\let\mymatrixII\mymatrix
\usepackage{mymatrix}
\usepackage{longtable}
\usepackage{amssymb}
\begin{document}

\setcounter{MaxMatrixCols}{14}

\newcounter{bnomer} \newcounter{snomer}
\newcounter{bsnomer}
\setcounter{bnomer}{0}
\renewcommand{\thesnomer}{\thebnomer.\arabic{snomer}}
\renewcommand{\thebsnomer}{\thebnomer.\arabic{bsnomer}}
\renewcommand{\refname}{\begin{center}\large{\textbf{References}}\end{center}}

\newcommand{\sect}[1]{%
\setcounter{snomer}{0}\setcounter{bsnomer}{0}
\refstepcounter{bnomer}
\par\bigskip\begin{center}\large{\textbf{\arabic{bnomer}. {#1}}}\end{center}}
\newcommand{\sst}{%
\refstepcounter{bsnomer}
\par\bigskip\textbf{\arabic{bnomer}.\arabic{bsnomer}. }}
\newcommand{\defi}[1]{%
\refstepcounter{snomer}
\par\medskip\textbf{Definition \arabic{bnomer}.\arabic{snomer}. }{#1}\par\medskip}
\newcommand{\theo}[2]{%
\refstepcounter{snomer}
\par\textbf{Теорема \arabic{bnomer}.\arabic{snomer}. }{#2} {\emph{#1}}\hspace{\fill}$\square$\par\medskip}
\newcommand{\mtheop}[2]{%
\refstepcounter{snomer}
\par\textbf{Theorem \arabic{bnomer}.\arabic{snomer}. }{\emph{#1}}
\par\textsc{Proof}. {#2}\hspace{\fill}$\square$\par}
\newcommand{\mcorop}[2]{%
\refstepcounter{snomer}
\par\textbf{Corollary \arabic{bnomer}.\arabic{snomer}. }{\emph{#1}}
\par\textsc{Proof}. {#2}\hspace{\fill}$\square$\par\medskip}
\newcommand{\mtheo}[1]{%
\refstepcounter{snomer}
\par\medskip\textbf{Theorem \arabic{bnomer}.\arabic{snomer}. }{\emph{#1}}\par\medskip}
\newcommand{\mlemm}[1]{%
\refstepcounter{snomer}
\par\medskip\textbf{Lemma \arabic{bnomer}.\arabic{snomer}. }{\emph{#1}}\par\medskip}
\newcommand{\mprop}[1]{%
\refstepcounter{snomer}
\par\medskip\textbf{Proposition \arabic{bnomer}.\arabic{snomer}. }{\emph{#1}}\par\medskip}
\newcommand{\theobp}[2]{%
\refstepcounter{snomer}
\par\medskip\textbf{Теорема \arabic{bnomer}.\arabic{snomer}. }{#2} {\emph{#1}}\par\medskip}
\newcommand{\theop}[2]{%
\refstepcounter{snomer}
\par\medskip\textbf{Theorem \arabic{bnomer}.\arabic{snomer}. }{\emph{#1}}
\par\textsc{Proof}. {#2}\hspace{\fill}$\square$\par\medskip}
\newcommand{\theosp}[2]{%
\refstepcounter{snomer}
\par\medskip\textbf{Теорема \arabic{bnomer}.\arabic{snomer}. }{\emph{#1}}
\par\textbf{Схема доказательства}. {#2}\hspace{\fill}$\square$\par\medskip}
\newcommand{\exam}[1]{%
\refstepcounter{snomer}
\par\medskip\textbf{Example \arabic{bnomer}.\arabic{snomer}. }{#1}\par\medskip}
\newcommand{\deno}[1]{%
\refstepcounter{snomer}
\par\medskip\textbf{Definition \arabic{bnomer}.\arabic{snomer}. }{#1}\par\medskip}
\newcommand{\post}[1]{%
\refstepcounter{snomer}
\par\medskip\textbf{Предложение \arabic{bnomer}.\arabic{snomer}. }{\emph{#1}}\hspace{\fill}$\square$\par}
\newcommand{\postp}[2]{%
\refstepcounter{snomer}
\par\medskip\medskip\textbf{Proposition \arabic{bnomer}.\arabic{snomer}. }{\emph{#1}}%
\ifhmode\par\fi\textsc{Proof}. {#2}\hspace{\fill}$\square$\par\medskip}
\newcommand{\lemm}[1]{%
\refstepcounter{snomer}
\par\medskip\textbf{Lemma \arabic{bnomer}.\arabic{snomer}. }{\emph{#1}}\hspace{\fill}$\square$\par\medskip}
\newcommand{\lemmp}[2]{%
\refstepcounter{snomer}
\par\medskip\textbf{Lemma \arabic{bnomer}.\arabic{snomer}. }{\emph{#1}}
\par\textsc{Proof}. {#2}\hspace{\fill}$\square$\par\medskip}
\newcommand{\coro}[1]{%
\refstepcounter{snomer}
\par\medskip\textbf{Corollary \arabic{bnomer}.\arabic{snomer}. }{\emph{#1}}\hspace{\fill}$\square$\par\medskip}
\newcommand{\mcoro}[1]{%
\refstepcounter{snomer}
\par\medskip\textbf{Corollary \arabic{bnomer}.\arabic{snomer}. }{\emph{#1}}\par\medskip}
\newcommand{\corop}[2]{%
\refstepcounter{snomer}
\par\medskip\textbf{Corollary \arabic{bnomer}.\arabic{snomer}. }{\emph{#1}}
\par\textsc{Proof}. {#2}\hspace{\fill}$\square$\par\medskip}
\newcommand{\nota}[1]{%
\refstepcounter{snomer}
\par\medskip\textbf{Remark \arabic{bnomer}.\arabic{snomer}. }{#1}\par\medskip}
\newcommand{\propp}[2]{%
\refstepcounter{snomer}
\par\medskip\textbf{Proposition \arabic{bnomer}.\arabic{snomer}. }{\emph{#1}}
\par\textsc{Proof}. {#2}\hspace{\fill}$\square$\par\medskip}
\newcommand{\hypo}[1]{%
\refstepcounter{snomer}
\par\medskip\textbf{Conjecture \arabic{bnomer}.\arabic{snomer}. }{\emph{#1}}\par\medskip}
\newcommand{\prop}[1]{%
\refstepcounter{snomer}
\par\medskip\textbf{Proposition \arabic{bnomer}.\arabic{snomer}. }{\emph{#1}}\hspace{\fill}$\square$\par\medskip}

\newcommand{\Ind}[3]{%
\mathrm{Ind}_{#1}^{#2}{#3}}
\newcommand{\Res}[3]{%
\mathrm{Res}_{#1}^{#2}{#3}}
\newcommand{\epsi}{\epsilon}
\newcommand{\tri}{\triangleleft}
\newcommand{\Supp}[1]{%
\mathrm{Supp}(#1)}

\newcommand{\reg}{\mathrm{reg}}
\newcommand{\sreg}{\mathrm{sreg}}
\newcommand{\codim}{\mathrm{codim}\,}
\newcommand{\chara}{\mathrm{char}\,}
\newcommand{\rk}{\mathrm{rk}\,}
\newcommand{\id}{\mathrm{id}}
\newcommand{\Ad}{\mathrm{Ad}}
\newcommand{\col}{\mathrm{col}}
\newcommand{\row}{\mathrm{row}}
\newcommand{\low}{\mathrm{low}}
\newcommand{\pho}{\hphantom{\quad}\vphantom{\mid}}
\newcommand{\fho}[1]{\vphantom{\mid}\setbox0\hbox{00}\hbox to \wd0{\hss\ensuremath{#1}\hss}}
\newcommand{\wt}{\widetilde}
\newcommand{\wh}{\widehat}
\newcommand{\ad}[1]{\mathrm{ad}_{#1}}
\newcommand{\tr}{\mathrm{tr}\,}
\newcommand{\GL}{\mathrm{GL}}
\newcommand{\SO}{\mathrm{SO}}
\newcommand{\Or}{\mathrm{O}}
\newcommand{\Sp}{\mathrm{Sp}}
\newcommand{\Mat}{\mathrm{Mat}}
\newcommand{\Stab}{\mathrm{Stab}}
\newcommand{\empr}[2]{[-{#1},{#1}]\times[-{#2},{#2}]}

\newcommand{\vfi}{\varphi}
\newcommand{\teta}{\vartheta}
\newcommand{\lee}{\leqslant}
\newcommand{\gee}{\geqslant}
\newcommand{\Fp}{\mathbb{F}}
\newcommand{\Rp}{\mathbb{R}}
\newcommand{\Zp}{\mathbb{Z}}
\newcommand{\Cp}{\mathbb{C}}
\newcommand{\ut}{\mathfrak{u}}
\newcommand{\at}{\mathfrak{a}}
\newcommand{\nt}{\mathfrak{n}}
\newcommand{\spt}{\mathfrak{sp}}
\newcommand{\htt}{\mathfrak{h}}
\newcommand{\slt}{\mathfrak{sl}}
\newcommand{\ot}{\mathfrak{so}}
\newcommand{\rt}{\mathfrak{r}}
\newcommand{\rad}{\mathfrak{rad}}
\newcommand{\bt}{\mathfrak{b}}
\newcommand{\gt}{\mathfrak{g}}
\newcommand{\vt}{\mathfrak{v}}
\newcommand{\pt}{\mathfrak{p}}
\newcommand{\Po}{\EuScript{P}}
\newcommand{\Uo}{\EuScript{U}}
\newcommand{\Fo}{\EuScript{F}}
\newcommand{\Do}{\EuScript{D}}
\newcommand{\Eo}{\EuScript{E}}
\newcommand{\So}{\EuScript{S}}
\newcommand{\Iu}{\mathcal{I}}
\newcommand{\Xu}{\mathcal{X}}
\newcommand{\Mo}{\mathcal{M}}
\newcommand{\Nu}{\mathcal{N}}
\newcommand{\Ro}{\mathcal{R}}
\newcommand{\Co}{\mathcal{C}}
\newcommand{\Lo}{\mathcal{L}}
\newcommand{\Ou}{\mathcal{O}}
\newcommand{\Au}{\mathcal{A}}
\newcommand{\Vu}{\mathcal{V}}
\newcommand{\Bu}{\mathcal{B}}
\newcommand{\Sy}{\mathcal{Z}}
\newcommand{\Sb}{\mathcal{F}}
\newcommand{\Gr}{\mathcal{G}}

\author{Mikhail V. Ignatyev\thanks{The work on Section~\ref{sect:proof_mtheo} was performed at the NRU HSE with the support from the Russian Science Foundation, grant no.\ 16--41--01013. The work on Section~\ref{sect:concluding_remarks} has been supported by RFBR grant no.\ 16--01--00154a and by\break ISF grant no. 797/14.}}
\date{\small Samara National Research University\\
443011, Ak. Pavlova, 1, Samara, Russia\\
\texttt{mihail.ignatev@gmail.com}}
\title{\Large{On involutions in the Weyl group and $B$-orbit closures in the orthogonal case}} \maketitle

\vspace{-0.5cm}\begin{center}
\begin{tabular}{p{15cm}}
\small{\textsc{Abstract}. We study coadjoint $B$-orbits on $\nt^*$, where $B$ is a Borel subgroup of a complex orthogonal group $G$, and $\nt$ is the Lie algebra of the unipotent radical of $B$. To each basis involution $w$ in the Weyl group $W$ of $G$ one can assign the associated $B$-orbit $\Omega_w$. We prove that, given basis involutions $\sigma$, $\tau$ in $W$, if the orbit $\Omega_{\sigma}$ is contained in the closure of the orbit $\Omega_{\tau}$ then $\sigma$ is less than or equal to $\tau$ with respect to the Bruhat order on $W$. For a basis involution $w$, we also compute the dimension of $\Omega_w$ and present a conjectural description of the closure of $\Omega_w$.}
\end{tabular}
\end{center}

\vspace{-0.7cm}\sect{Introduction, definitions and the main result}

\vspace{-0.5cm}\sst Let $G$ be a complex reductive algebraic group, $B$ a Borel subgroup of $G$, $\Phi$ the root system of~$G$ and $W=W(\Phi)$ the Weyl group of $\Phi$. It is well-known that the Bruhat order on~$W$ encodes the cell decomposition of the flag variety $G/B$ (see, e.g., \cite{BileyLakshmibai1}). Denote by $\Iu(\Phi)$ the poset of involutions in~$W$ (i.e., elements of~$W$ of order 2). In \cite{RichardsonSpringer1}, R. Richardson and T.~Springer showed that $\Iu(A_{2n})$ encodes the incidences among the closed $B$-orbits on the symmetric variety $\mathrm{SL}_{2n+1}(\Cp)/\mathrm{SO}_{2n+1}(\Cp)$. In~\cite{BagnoCherniavsky1}, E.~Bagno and Y.~{Cher\-nav\-sky} presented a~{geo\-met\-ri\-cal} {inter\-pre\-ta\-tion} of the poset $\Iu(A_n)$, considering the action of the Borel subgroup of $\GL_n(\Cp)$ on symmetric matrices by congruence. F.~Incitti studied the poset $\Iu(\Phi)$ from a purely combinatorial point of view for the case of classical root system~$\Phi$ (see \cite{Incitti1}, \cite{Incitti2}, \cite{Incitti3}). In particular, he proved that this poset is graded, calculated the rank function and described the covering relation.

In \cite{Ignatev3}, we presented another geometrical interpretation of $\Iu(A_{n-1})$ in terms of \emph{coadjoint} $B$-orbits. Precisely, let $U$ be the the unipotent radical of $B$, and Let $\nt$ be the Lie algebra of $U$. Since $B$ acts on $\nt$ via the adjoint action, one can consider the dual action of $B$ on $\nt^*$, which is called \emph{coadjoint}. To each involution $\sigma\in\Iu(A_{n-1})$ one can assign the $B$-orbit $\Omega_{\sigma}\subseteq\nt^*$ (see Subsections~\ref{sst:A_n_C_n_definitions},~\ref{sst:B_n_D_n_definitions} for precise definitions). By~\cite[Theorem 1.1]{Ignatev3}, for $G=\GL_n(\Cp)$, $\Omega_{\sigma}$ is contained in the Zariski closure $\overline{\Omega}_{\tau}$ of $\Omega_{\tau}$ if and only if $\sigma$ is less or equal to $\tau$ with respect to the Bruhat order $\leq_B$. In \cite{Ignatyev4}, completely similar results were obtained for $G=\Sp_{2n}(\Cp)$ (i.e., for $\Iu(C_n)$), see \cite[Theorem 1.1]{Ignatyev4}. In some sense, these results are ``dual'' to A. Melnikov's results \cite{Melnikov1}, \cite{Melnikov2}, \cite{Melnikov3}.

In this paper, we establish similar results for the cases $\Phi=B_n$ and $\Phi=D_n$. Namely, let $G$ be the orthogonal group of rank $n$, i.e., $G=\Or_{2n+1}(\Cp)$ or $\Or_{2n}(\Cp)$ (respectively, $\Phi=B_n$ or $D_n$). In general, $\Omega_{\sigma}\subseteq\overline{\Omega}_{\tau}$ is \emph{not} equivalent to $\sigma\leq_B\tau$, see Examples~\ref{exam:if_contained_not_Bruhat_B_n} and \ref{exam:if_Bruhat_not_contained_B_n} below. On the other hand, we believe that these conditions \emph{are} equivalent if we restrict ourselves to the case of so-called basis involutions. An involution~$\sigma$ is called \emph{basis} if there are no $i$ such that $\sigma(i)=-i$. (Here $W$ is standardly identified with certain subgroup of the symmetric group $S_{\pm n}$ on the $2n$ letters $1,~\ldots,~n,~-n,~\ldots,~-1$, see Subsection~\ref{sst:A_n_C_n_definitions} for the precise definition). The main result of the paper is as follows.\newpage

\mtheo{Let $\sigma$, $\tau$ be basis involutions in the Weyl group $W$ of type $B_n$ or $D_n$. If the orbit~$\Omega_{\sigma}$ is contained in the Zariski closure of the orbit $\Omega_{\tau}$, then $\sigma$ is less or equal to $\tau$ with respect to the Bruhat order on the group $W$.\label{mtheo}}

The paper is organized as follows. In the rest of this section we briefly recall basic facts about classical groups and collect our previous results about $\Iu(A_{n-1})$ and $\Iu(C_n)$, see Subsection \ref{sst:A_n_C_n_definitions}. Then, we give precise definitions for the orthogonal case and formulate some Incitti's results about involutions needed in the sequel, see Subsection~\ref{sst:B_n_D_n_definitions}.

Section~\ref{sect:proof_mtheo} is devoted to the proof of Theorem~\ref{mtheo}. Precisely, in Subsection~\ref{sst:if_contained_then_Bruhat_B_n} we prove it for $B_n$, see Theorem~\ref{theo:if_contained_then_Bruhat_B_n}. Next, in Subsection~\ref{sst:if_contained_then_Bruhat_D_n} we prove this theorem for $D_n$ (this requires some additional work due to the fact that the Bruhat order in this case has more complicated description than for~$B_n$). In Subsection~\ref{sst:if_Bruhat_then_contained} we discuss the equivalence of the conditions $\Omega_{\sigma}\subseteq\overline{\Omega}_{\tau}$ and $\sigma\leq_B\tau$ for basis involutions. Namely, using Incitti's results, we present a conjectural way how to prove that if $\sigma\leq_B\tau$ then $\Omega_{\sigma}$ is contained in the closure of ${\Omega}_{\tau}$.

Finally, in Section~\ref{sect:concluding_remarks} we discuss some related facts and conjectures. In Subsection \ref{sst:dim_Omega}, we obtain a~formula for the dimension of the orbit $\Omega$ (see Theorem~\ref{theo:dim_Omega}). In Subsection \ref{sst:Schubert}, a conjectural approach to orbits associated with involutions in terms of tangent cones to Schubert subvarieties of the flag variety $G/B$ is presented.

\medskip\textsc{Acknowledgements}. A part of this work (Section~\ref{sect:concluding_remarks}) was done during my stay at University of Haifa. I~would like to express my gratitude to Prof. Dr. Anna Melnikov for her hospitality and fruitful discussions.

\sst\label{sst:A_n_C_n_definitions} From now on and to the end of the paper $G$ denotes one of the classical complex algebraic groups $\GL_n(\Cp)$, $\Or_{2n+1}(\Cp)$, $\Sp_{2n}(\Cp)$ or $\Or_{2n}(\Cp)$. The group $\Or_{2n+1}(\Cp)$ (respectively, $\Sp_{2n}(\Cp)$, $\Or_{2n}(\Cp)$) is realized as the subgroup of $\GL_{2n+1}(\Cp)$ (respectively, of $\GL_{2n}(\Cp)$) consisting of all invertible matrices~$g$ such that $$\beta(gu,gv)=\beta(u,v)$$ for all $u,v$ in $\Cp^{2n+1}$ (respectively, in $\Cp^{2n}$), where $\beta$ is the bilinear form on $\Cp^{2n+1}$ (respectively, on $\Cp^{2n}$) defined as follows:
\begin{equation*}
\beta(u,v)=\begin{cases}
u_0v_0+\sum\limits_{i=1}^n(u_iv_{-i}+u_{-i}v_i)&\text{for }\Or_{2n+1}(\Cp),\\
\sum\limits_{i=1}^n(u_iv_{-i}-u_{-i}v_i)&\text{for }\Sp_{2n}(\Cp),\\
\sum\limits_{i=1}^n(u_iv_{-i}+u_{-i}v_i)&\text{for }\Or_{2n}(\Cp).
\end{cases}
\end{equation*}
Here for $\Or_{2n+1}(\Cp)$ (respectively, for $\Sp_{2n}(\Cp)$ and $\Or_{2n}(\Cp)$) we denote by $e_1,\ldots,e_n,e_{-n},\ldots,e_{-1}$ (res\-pec\-tively, by $e_1,\ldots,e_n,e_0,e_{-n},\ldots,e_{-1}$ and $e_1,\ldots,e_n,e_{-n},\ldots,e_{-1}$) the standard basis of $\Cp^{2n+1}$ (res\-pec\-tively, of $\Cp^{2n}$), and by $x_i$ the coordinate of a vector $x$ corresponding to $e_i$.

The set of all diagonal matrices from $G$ is a maximal torus in $G$; we denote it by $H$. Let $\Phi$ be the root system of $G$ with respect to $H$. Note that $\Phi$ is of type $A_{n-1}$ (respectively, $B_n$, $C_n$ and $D_n$) for $\GL_n(\Cp)$ (respectively, for $\Or_{2n+1}(\Cp)$, $\Sp_{2n}(\Cp)$ and $\Or_{2n}(\Cp)$). The set of all upper-triangular matrices from $G$ is a Borel subgroup of $G$ containing $H$; we denote it by $B$. Let $\Phi^+$ be the set of positive roots with respect to $B$. As usual, we identify $\Phi^+$ with the following subset of the $n$-dimensional Euclidean space $\Rp^n$ with the standard inner product (see, e.g., \cite{Bourbaki1}):
\begin{equation*}\predisplaypenalty=0
\begin{split}
A_{n-1}^+&=\{\epsi_i-\epsi_j,~1\leq i<j\leq n\},\\
B_n^+&=\{\epsi_i-\epsi_j,~1\leq i<j\leq n\}\cup\{\epsi_i+\epsi_j,~1\leq i<j\leq n\}\cup\{\epsi_i,~1\leq i\leq n\},\\
C_n^+&=\{\epsi_i-\epsi_j,~1\leq i<j\leq n\}\cup\{\epsi_i+\epsi_j,~1\leq i<j\leq n\}\cup\{2\epsi_i,~1\leq i\leq n\},\\
D_n^+&=\{\epsi_i-\epsi_j,~1\leq i<j\leq n\}\cup\{\epsi_i+\epsi_j,~1\leq i<j\leq n\}.\\
\end{split}
\end{equation*}
Here $\{\epsi_i\}_{i=1}^n$ is the standard basis of $\Rp^n$.\newpage

Denote by $U$ the group of all strictly upper-triangular matrices from $G$ with $1$'s on the diagonal, then $U$ is the unipotent radical of $B$. Let $\gt$, $\htt$, $\bt$, $\nt$ be the Lie algebras of $G$, $H$, $B$, $U$ respectively. Then $\nt$ has a basis consisting of root vectors $e_{\alpha}$, $\alpha\in\Phi^+$, where
\begin{equation*}\predisplaypenalty=10000
\begin{split}
e_{\epsi_i}&=\sqrt{2}(e_{0,i}-e_{-i,0}),~e_{2\epsi_i}=e_{i,-i},\\
e_{\epsi_i-\epsi_j}&=\begin{cases}
e_{i,j}&\text{for }A_{n-1},\\
e_{i,j}-e_{-j,-i}&\text{for }B_n,~C_n\text{ and }D_n,
\end{cases}\\
e_{\epsi_i+\epsi_j}&=\begin{cases}
e_{i,-j}-e_{j,-i}&\text{for }B_n\text{ and }D_n,\\
e_{i,-j}+e_{j,-i}&\text{for }C_n,
\end{cases}
\end{split}
\end{equation*}
and $e_{i,j}$ are the usual elementary matrices. For $\Or_{2n+1}(\Cp)$ (respectively, for $\Sp_{2n}(\Cp)$ and $\Or_{2n}(\Cp)$) we index the rows (from left to right) and the columns (from top to bottom) of matrices by the numbers $1,~\ldots~,~n,~0,~-n,~\ldots,~-1$ (respectively, by the numbers $1,~\ldots,~n,~-n,~\ldots,~-1$). Note that $\gt=\htt\oplus\nt\oplus\nt_-$, where $\nt_-=\langle e_{-\alpha},~\alpha\in\Phi^+\rangle_{\Cp}$, and, by definition, $e_{-\alpha}=e_{\alpha}^T$. (The superscript~$T$ always indicates matrix transposition.)

Since $\{e_{\alpha},~\alpha\in\Phi^*\}$ is a basis of $\nt$, one can consider the dual basis $\{e_{\alpha}^*,~\alpha\in\Phi^+\}$ of the dual space~$\nt^*$. The group $B$ acts on $\nt$ by the adjoint action (actually, by conjugation), so there exists the dual (\emph{coadjoint}) action of $B$ on $\nt^*$. We will denote the result of this action by $g.\lambda$ for $g\in B$, $\lambda\in\nt^*$. By definition,
\begin{equation*}
\langle g.\lambda,x\rangle=\langle\lambda,g^{-1}xg\rangle,\text{ }g\in B,\text{ }x\in\nt,\text{ }\lambda\in\nt^*.
\end{equation*}
Orbits of the coadjoint action play the crucial role in representation theory of the groups $B$ and $U$, see, e.g., \cite{Kirillov1}, \cite{Kirillov2}.

It is very convenient to identify $\nt^*$ with the space $\nt_-$ via $$\langle\lambda,x\rangle=\tr\lambda x,\text{ }\lambda\in\nt_-,\text{ }x\in\nt.$$
For this reason, we will denote $\nt_-$ by $\nt^*$ and interpret it as the dual space of $\nt$. Under this identification,
\begin{equation*}
e_{\alpha}^*=\begin{cases}\label{formula:Killing}
e_{\alpha}^T&\text{if }\Phi=A_{n-1},\text{ or }\Phi=C_n\text{ and }\alpha=2\epsi_i,\\
e_{\alpha}^T/4&\text{if }\Phi=B_n\text{ and }\alpha=\epsi_i,\\
e_{\alpha}^T/2&\text{otherwise}.
\end{cases}
\end{equation*}
Note that if $g\in B$, $\lambda\in\nt^*$, then $$g.\lambda=(g\lambda g^{-1})_{\low},$$
where $A_{\low}$ denotes the strictly lower-triangular part of a matrix $A$, i.e.,
\begin{equation*}
(A_{\low})_{i,j}=\begin{cases}A_{i,j}&\text{for }i>j,\\
0&\text{for }i\leq j.
\end{cases}
\end{equation*}

For a given $\lambda\in\nt^*$, let $\Omega_{\lambda}$ and $\Theta_{\lambda}$ denote its $B$-orbit and $U$-orbit under the coadjoint action respectively. A subset $D\subset\Phi^+$ is called \emph{orthogonal} if it consists of pairwise orthogonal roots. To each orthogonal subset $D$ and each map $\xi\colon D\to\Cp^{\times}$ one can assign the linear forms $$f_D=\sum_{\alpha\in D}e_{\alpha}^*\in\nt^*,~f_{D,\xi}=\sum_{\alpha\in D}\xi(\alpha)e_{\alpha}^*\in\nt^*.$$
(Obviously, $f_D=f_{D,\xi_1}$, where $\xi_1(\alpha)=1$ for all $\alpha\in D$.) Given an orthogonal subset $D\subseteq\Phi^+$, we say that the orbits $\Omega_D=\Omega_{f_D}$ and $\Theta_{D,\xi}=\Theta_{f_{D,\xi}}$ are \emph{associated} with $D$. Note that $U$-orbits associated with orthogonal subsets and their generalizations were studied, in particular, in \cite{Panov}, \cite{Ignatev1}, \cite{Ignatev2}, \cite{IgnatyevVasyukhin}.

Now, let $W$ be the Weyl group of $\Phi$. We denote by $s_{\alpha}$ the reflection in $W$ corresponding to a root~$\alpha$, and say that $s_{\alpha}$ is a \emph{simple reflection} if $\alpha$ is a simple root. For $\Phi=A_{n-1}$, $W\cong S_n$ is isomorphic to the symmetric group on the $n$ letters $1,~\ldots,~n$ via the isomorphism $s_{\epsi_i-\epsi_j}\mapsto(i,j)$, where $(i,j)$ is the transposition interchanging $i$ and $j$. For other classical root systems, denote by $S_{\pm n}$ the symmetric group on the $2n$ letters $1,~\ldots,~n,~-n,~\ldots,~-1$ and consider the monomorphism from $W$ to $S_{\pm n}$ defined by the formulas
\begin{equation*}
\begin{split}
&s_{\epsi_i-\epsi_j}\mapsto(i,j)(-j,-i),\\
&s_{\epsi_i+\epsi_j}\mapsto(i,-j)(j,-i),\\
&s_{\epsi_i}\mapsto(i,-i),~s_{2\epsi_i}\mapsto(i,-i).
\end{split}
\end{equation*}
For $B_n$ and $C_n$, the image of this monomorphism coincides with the \emph{hyperoctahedral} group, that is, the subgroup of $S_{\pm n}$  consisting of all permutations $w$ from $S_{\pm n}$ such that $w(-i)=-w(i)$ for each $1\leq i\leq n$. For $D_n$, the image of this monomorphism coincides with the \emph{even-signed hyperoctahedral} group, that is, the subgroup of $S_{\pm n}$ consisting of all $w\in S_{\pm n}$ such that $w(-i)=-w(i)$ for each $1\leq i\leq n$ and the number ${|\{i>0\mid w(i)<0\}|}$ is even. We will identify $W$ with its image under the above monomorphism.

\nota{i) Note that every $w\in W$ is completely determined by its restriction to the subset $\{1,\ldots,n\}$. This allows us to use the usual two-line notation: if $w(i)=w_i$ for $1\leq i\leq n$, then we will write $w=\begin{pmatrix}1&2&\ldots&n\\w_1&w_2&\ldots&w_n\end{pmatrix}$. For instance, if $\Phi=D_5$, then $$s_{\epsi_1+\epsi_5}s_{\epsi_2+\epsi_4}s_{\epsi_2-\epsi_4}=\begin{pmatrix}1&2&3&4&5\\-5&-2&3&-4&-1\end{pmatrix}.$$

ii) Note also that the set of simple roots has the following form: $\Delta=\{\alpha_1,\ldots,\alpha_n\}$, where\break $\alpha_1=\epsi_1-\epsi_2$, $\ldots$, $\alpha_{n-1}=\epsi_{n-1}-\epsi_n$, and
\begin{equation*}
\alpha_n=\begin{cases}
\epsi_n&\text{for }B_n,\\
2\epsi_n&\text{for }C_n,\\
\epsi_{n-1}+\epsi_{n+1}&\text{for }D_n.
\end{cases}
\end{equation*}}

Recall that a \emph{reduced decomposition} of an element $w\in W$ is an expression of $w$ as a product of simple reflections of minimal possible length. Given $v,~w\in W$, we say that $v$ is less or equal to $w$ with respect to the \emph{Bruhat order}, written $v\leq_B w$, if some reduced decomposition for $v$ is a subword of some reduced decomposition for $w$. It is well-known that this order plays the crucial role in many geometric aspects of theory of algebraic groups. For instance, the Bruhat order encodes the incidences among Schubert varieties

From now on and to the end of this subsection, let $G=\GL_n(\Cp)$ or $\Sp_{2n}(\Cp)$, i.e., $\Phi=A_{n-1}$ or $C_n$ respectively. There exists a nice combinatorial description of the Bruhat order on $W$. First, consider the case $A_{n-1}$. Given $w\in W$, denote by $X_w$ the $n\times n$ matrix defined by
\begin{equation*}
(X_w)_{i,j}=\begin{cases}1,&\text{ if }w(j)=i,\\
0&\text{ otherwise}.
\end{cases}
\end{equation*}
It is called the $0$--1 matrix, permutation matrix or \emph{rook placement} for $w$. Define the matrix $R_w$ by putting its $(i,j)$th element to be equal to the rank of the lower left $(n-i+1)\times j$ submatrix of $X_w$. In other words, $(R_w)_{i,j}$ is just the number or rooks located non-strictly to the South-West from $(i,j)$.\newpage

\exam{Let $n=6$, $w=s_{\epsi_1-\epsi_4}s_{\epsi_3-\epsi_5}=\begin{pmatrix}1&2&3&4&5&6\\4&2&5&1&3&6\end{pmatrix}$. Here we draw the matrices $X_w$ and $R_w$ (rooks are marked by $\otimes$):
\begin{equation*}X_w=
\mymatrix{
\pho& \pho& \pho& \otimes& \pho& \pho\\
\pho& \otimes& \pho& \pho& \pho& \pho\\
\pho& \pho& \pho& \pho& \otimes& \pho\\
\otimes& \pho& \pho& \pho& \pho& \pho\\
\pho& \pho& \otimes& \pho& \pho& \pho\\
\pho& \pho& \pho& \pho& \pho& \otimes\\
}\ ,\ R_w=
\mymatrix{
\fho1& \fho2& \fho3& \fho4& \fho5& \fho6\\
\fho1& \fho2& \fho2& \fho3& \fho4& \fho5\\
\fho1& \fho1& \fho2& \fho2& \fho3& \fho4\\
\fho1& \fho1& \fho2& \fho2& \fho2& \fho3\\
\fho0& \fho0& \fho1& \fho1& \fho1& \fho2\\
\fho0& \fho0& \fho0& \fho0& \fho0& \fho1\\
}\ .
\end{equation*}\label{exam:Bruhat_A_n}}

Given two arbitrary matrices $X$, $Y$ of the same size with integer entries, we will write $X\leq Y$ if $X_{i,j}\leq Y_{i,j}$ for all $i$, $j$. It turns out that, for $\Phi=A_{n-1}$ and $v$, $w\in W$, $v\leq_Bw$ if and only if $R_v\leq R_w$, see, e.g., \cite[Theorem 2.1.5]{BjornerBrenti}.

For $C_n$, the description of the Bruhat order is very similar. Precisely, to each $w\in W$ one can assign the matrices $X_w$ and $R_w$ exactly by the same rule. (Here $X_w$ and $R_w$ are $2n\times2n$ matrices, which rows and columns are indexed by the numbers $1,~\ldots,~n,~-n,~\ldots,~-1$.) According to\break \cite[Theorem 8.1.1]{BjornerBrenti}, for $\Phi=C_n$ and $v,~w\in W$, $v\leq_Bw$ if and only if $R_v\leq R_w$, as above. In other words, the Bruhat order on $W$ is nothing but the restriction of the Bruhat order on $S_{\pm n}$ to $W$.

Let $G=\GL_n(\Cp)$ or $\Sp_{2n}(\Cp)$, and $w$ be an involution in $W$, i.e., $w\in\Iu(\Phi)$. Then $w$ can be expressed as a product of pairwise commuting reflections. In other words, there exists an orthogonal subset $D\subseteq\Phi^+$ such that $w=\prod_{\alpha\in D}s_{\alpha}$. For $A_{n-1}$, such an expression is clearly unique. For $C_n$, such an expression is unique if we require $D$ to be strongly orthogonal, which means that, given $\alpha,~\beta\in D$, neither $\alpha+\beta\in\Phi^+$ nor $\alpha-\beta\in\Phi^+$. Thus, in both cases, the subset $D$ is uniquely determined by $w$. We call the subset $D$ the \emph{support} of the involution $w$ and denote it by $D=\Supp{w}$. For instance, if $\Phi=C_6$ and
\begin{equation*}
w=\begin{pmatrix}1&2&3&4&5&6\\3&-6&1&-4&-5&-2\end{pmatrix}=s_{\epsi_1-\epsi_3}s_{\epsi_2+\epsi_6}s_{2\epsi_4}s_{2\epsi_5}=
s_{\epsi_1-\epsi_3}s_{\epsi_2+\epsi_6}s_{\epsi_4-\epsi_5}s_{\epsi_4+\epsi_5},
\end{equation*}
then $\Supp{w}=\{\epsi_1-\epsi_3,~\epsi_2+\epsi_6,~2\epsi_4,~2\epsi_5\}$.

\defi{Let\label{defi:orbits_associated} $w\in\Iu(\Phi)$, and $\xi\colon D\to\Cp^{\times}$ be a map, where $D=\Supp{w}$. Denote $f_w=f_D$, $f_{w,\xi}=f_{D,\xi}$, $\Omega_w=\Omega_D$ and $\Theta_{w,\xi}=\Theta_{D,\xi}$. We say that the orbits $\Omega_w$ and $\Theta_{w,\xi}$ are associated with the involution~$w$. Note that, thanks to \cite[Lemma 2.1]{Ignatev3} and \cite[Lemma 1.8]{Ignatyev4}, $$\Omega_w=\bigcup_{\xi}\Theta_{w,\xi}.$$}

To formulate the description of the incidences among $B$-orbits associated with involutions, we need one more partial order on $\Iu(\Phi)$. Namely, given $w\in W$, we put $R_w^*=(R_w)_{\low}$ and write $v\leq^*w$ for $v,~w\in\Iu(\Phi)$ if $R_v^*\leq R_w^*$. Then, according to \cite[Theorem 1.1]{Ignatev3} and \cite[Theorem 1.1]{Ignatyev4}, we have the following result.

\mtheo{Let\label{mtheo:A_n_C_n} $\Phi=A_{n-1}$ or $C_n$\textup, and $v,~w\in\Iu(\Phi)$. Then the following conditions are equivalent\textup:
\begin{equation*}
\begin{split}
&\text{\textup{i) }}\Omega_v\subseteq\overline{\Omega}_w;\\
&\text{\textup{ii) }}v\leq^*w;\\
&\text{\textup{iii) }}v\leq_Bw.
\end{split}
\end{equation*}}

\sst\label{sst:B_n_D_n_definitions} Suppose now that $G=\Or_{2n+1}(\Cp)$ or $\Or_{2n(\Cp)}$, i.e., $\Phi=B_n$ or $D_n$ respectively. Since the Weyl group of $B_n$ is isomorphic to the Weyl group of $C_n$, the Bruhat order on $W$ for $B_n$ can be described completely similarly to the case $\Phi=C_n$: given $v,~w\in W$, one has $v\leq_Bw$ if and only if $R_v\leq R_w$, where $R_v$, $R_w$ are the $2n\times 2n$ defined in the previous subsection.

For $D_n$, the description of the Bruhat order is quite more complicated. Let $w\in W$. Given $a,b\in\{1,2,\ldots,n\}$, we say that $\empr{a}{b}$ is an \emph{empty rectangle} for $w$, if $$\{i\in[\pm n]\mid |i|\geq b\text{ and }|w(i)|\geq a\}=\varnothing.$$ Here $[\pm n]=\{1,\ldots,n,-n,\ldots,-1\}$. For instance, let $n=4$, $w=\begin{pmatrix}1&2&3&4\\-2&4&1&-3\end{pmatrix}$, then
\begin{equation*}
X_w=
\mymatrixII{
\pho& \pho& \otimes& \pho& \pho& \pho& \pho& \pho\\
\pho& \pho& \pho& \pho& \pho& \pho& \pho& \otimes\\
\pho& \pho& \pho& \pho& \otimes& \pho& \pho& \pho\\
\pho& \otimes& \pho& \pho& \pho& \pho& \pho& \pho\\
\pho& \pho& \pho& \pho& \pho& \pho& \otimes& \pho\\
\pho& \pho& \pho& \otimes& \pho& \pho& \pho& \pho\\
\otimes& \pho& \pho& \pho& \pho& \pho& \pho& \pho\\
\pho& \pho& \pho& \pho& \pho& \otimes& \pho& \pho\\}\ ,
\end{equation*}
and $\empr{4}{3}$, $\empr{4}{4}$ are empty rectangles for $w$. It turns out that given $v$, $w\in W$, $v\leq w$ if and only if
\begin{equation}
\begin{split}
&\text{i) }R_v\leq R_w;\\
&\text{ii) for all $a,b\in\{1,\ldots,n\}$, if $\empr{a}{b}$ is an empty rectangle}\\
&\hphantom{\text{ii) }}\text{for both $v$ and $w$ and $(R_v)_{-(a-1),b-1}=(R_w)_{-(a-1),b-1}$,}\\
&\hphantom{\text{ii) }}\text{then $(R_v)_{-(a-1),n}\equiv(R_w)_{-(a-1),n}\pmod{2}$.}\\
\end{split}\label{formula:Bruhat_Dn}
\end{equation}(See, e.g., \cite[Theorem 8.2.8]{BjornerBrenti}.)

\defi{Let $\Phi=B_n$ or $D_n$, and $w\in\Iu(\Phi)$. The involution $w$ is called \emph{basis} if $$|\{i\in\{1,~\ldots,~n\}\mid w(i)=-i\}|=0.$$}

We will denote the set of all basis involutions in $W$ by $\Bu(\Phi)$. It is clear that if $w$ is a basis involution then there exists the unique orthogonal subset $D\subseteq\Phi^+$ such that $$w=\prod_{\alpha\in D}s_{\alpha}.$$ As above, we call~$D$ \emph{the support} of $w$ and denote $D=\Supp{w}$. For example, for $\Phi=B_5$, $$w=\begin{pmatrix}1&2&3&4&5\\-5&2&4&3&-1\end{pmatrix},$$ we have $\Supp{w}=\{\epsi_1+\epsi_5,~\epsi_3-\epsi_4\}$. Now, given $w\in\Bu(\Phi)$ and a map $\xi\colon D=\Supp{w}\to\Cp^{\times}$, we define the linear forms $f_w$, $f_{w,\xi}$ and the orbits $\Omega_w$ and $\Theta_{w,\xi}$ exactly as in Definition~\ref{defi:orbits_associated}. As for $A_n$ and $C_n$, we say that $\Omega_w$ and $\Theta_{w,\xi}$ are \emph{associated} with $w$. The goal of the paper is to prove that, given $\sigma,~\tau\in\Bu(\Phi)$, $\Omega_{\sigma}\subseteq\overline{\Omega}_{\tau}$ implies $\sigma\leq_B\tau$. Note that this is clearly not true for arbitrary involutions in $W$, as it is shown in Examples~\ref{exam:if_contained_not_Bruhat_B_n}, \ref{exam:if_Bruhat_not_contained_B_n} below.

\sect{Proof of the main theorem}\label{sect:proof_mtheo}

\sst In\label{sst:if_contained_then_Bruhat_B_n} this subsection, we will check that, given $\sigma,~\tau\in\Bu(B_n)$, $\Omega_{\sigma}\subseteq\overline{\Omega}_{\tau}$ implies that $\sigma\leq_B\tau$.\break To do this, we need the following simple observation (cf. \cite[Subsection 3.3]{Melnikov1}, \cite[Lemma 2.1]{Ignatev3},\break \cite[Lemma 1.8]{Ignatyev4}, \cite[Lemma 1.3]{IgnatyevVasyukhin}).

\lemmp{Let\label{lemm:Omega_union_Theta} $\Phi=B_n$ or $D_n$, $w\in\Bu(\Phi)$ and $D=\Supp{w}$. Then $\Omega_w=\bigcup\Theta_{w,\xi}$, where the union is taken over all maps $\xi\colon D\to\Cp^{\times}$.}{It is well-known that the map $$\exp\colon\nt\to U,~x\mapsto\sum_{k\geq0}\dfrac{x^k}{k!}$$ is well-defined and is an isomorphism of affine varieties. For given $\alpha\in\Phi^+$, $s\in\Cp^{\times}$, put
\begin{equation*}\predisplaypenalty=0
\begin{split}
&x_{\alpha}(s)=\exp(se_{\alpha}),~x_{-\alpha}(s)=x_{\alpha}(s)^T,\\
&w_{\alpha}(s)=x_{\alpha}(s)x_{-\alpha}(-s^{-1})x_{\alpha}(s),\text{ }h_{\alpha}(s)=w_{\alpha}(s)w_{\alpha}(1)^{-1}.
\end{split}
\end{equation*}
Then $h_{\alpha}(s)$ is a diagonal matrix from $H$. Furthermore, the group $H$ is generated by $h_{\alpha}(s)$, $\alpha\in\Phi^+$, $s\in\Cp$, and $B=U\rtimes H$.

Let $\xi\colon D\to\Cp^{\times}$ be a map. To check that $\Theta_{w,\xi}\subseteq\Omega_w$, it is enough to find $h\in H$  such that $h.f_{w,\xi}=f_w$. One can easily see that if $\alpha\in D$ then $$h_{\alpha}(t).f_{w,\xi}=\sum_{\beta\in D,~\beta\neq\alpha}\xi(\beta)e_{\beta}^*+t^{-2}\xi(\alpha)e_{\alpha}^*.$$ Thus, $h.f_{w,\xi}=f_w$, where $h=\prod_{\alpha\in D}h_{\alpha}(\sqrt{\xi(\alpha)})$. (Here, given $s\in\Cp$, we denote by $\sqrt{s}$ a complex number such that $(\sqrt{s})^2=s$.)

On the other hand, let $h\in H$. We claim that $h.f_{w,\xi}=f_{w,\xi'}$ for some $\xi'$. Indeed, since $H$ is generated by $h_{\alpha}(s)$'s, $\alpha\in\Phi^+$, $s\in\Cp^{\times}$, we can assume without loss of generality that $h=h_{\alpha}(s)$ for some $\alpha$ and $s$. But in this case the statement follows immediately from the above. Since the group $B$ is isomorphic as an algebraic group to the semi-direct product $U\rtimes H$, for a~given $g\in B$, there exist unique $u\in U$, $h\in H$ such that $g=uh$. If $\xi\colon D\to\Cp^{\times}$ is the map such that $h.f_w=f_{w,\xi}$, then $g.f_w=u.f_{w,\xi}\in\Theta_{w,\xi}$. This concludes the proof.}

Now, if $\Phi=B_n$ or $D_n$, and $x\in\gt$, then, given $i,~j\in[\pm n]$, denote
\begin{equation*}
\pi_{i,j}(x)=\begin{pmatrix}x_{i,1}&\ldots&x_{i,j}\\
\vdots&\ddots&\vdots\\
x_{-1,1}&\ldots&x_{-1,j}\\
\end{pmatrix}.
\end{equation*}
It is easy to see that if $w\in\Bu(\Phi)$ then $\rk\pi_{i,j}(\lambda)=(R_w^*)_{i,j}$ for all $\lambda\in\Omega_w$ and all lower-triangular entries $(i,j)$ of $\lambda$ (cf. \cite[Lemma 2.2]{Ignatev3}, \cite[Lemma 2.4]{Ignatyev4}, \cite[Theorem 1.5]{IgnatyevVasyukhin}). Indeed, by definition of~$R_w^*$, $(R_w^*)_{i,j}=\rk\pi_{i,j}(f_w)$. Let $\xi\colon D\to\Cp^{\times}$ be a map. Since $$\rk\pi_{i,j}(f_{w,\xi})=\rk\pi_{i,j}(f_w)=(R_w^*)_{i,j},$$ it suffice to check that $\rk\pi_{i,j}(\lambda)=\rk\pi_{i,j}(u.\lambda)$ for $u\in U$, $\lambda\in\nt^*$. This follows immediately from the proof of \cite[Lemma 2.2]{Ignatev3}, because $u$ is an upper-triangular matrix with $1$'s on the diagonal and $\lambda$ is a~lower-triangular matrix with zeroes on the diagonal. Now we are ready to prove the main result of this subsection, cf. \cite[Proposition 2.3]{Ignatev3}, \cite[Proposition~2.5]{Ignatyev4}, \cite[Theorem 1.5]{IgnatyevVasyukhin}.\newpage

\theop{Let $\sigma,~\tau\in\Bu(B_n)$. If\label{theo:if_contained_then_Bruhat_B_n} $\Omega_{\sigma}\subseteq\overline{\Omega}_{\tau}$\textup, then $\sigma\leq_B\tau$.}{Suppose that $\sigma\nleq_B\tau$. According to Theorem~\ref{mtheo:A_n_C_n}, this means that there exist $i,~j$ such that $(R_{\sigma}^*)_{i,j}>(R_{\tau}^*)_{i,j}$. Denote
\begin{equation*}
Z=\{f\in\nt^*\mid\rk\pi_{r,s}(f)\leq(R_{\tau}^*)_{r,s}\text{ for
all }r,~s\}.
\end{equation*}
Clearly, $Z$ is closed with respect to the Zariski topology. It follows from the above that $\Omega_{\tau}\subseteq
Z$, so $\overline{\Omega}_{\tau}\subseteq Z$. But
$f_{\sigma}\notin Z$, hence $\Omega_{\sigma}\nsubseteq Z$, a
contradiction.}

\exam{If\label{exam:if_contained_not_Bruhat_B_n} involution $\sigma$ is not basis then, in general, its support (and, consequently, the associated $B$-orbit $\Omega_{\sigma}$) is not well-defined, because, in general, there are several different ways to represent $\sigma$ as a product of pairwise commuting reflections. For example, let $n=4$, then $$\sigma=\begin{pmatrix}1&2&3&4\\1&-2&-3&4\end{pmatrix}=s_{\epsi_2}s_{\epsi_3}=s_{\epsi_2-\epsi_3}s_{\epsi_2+\epsi_3}.$$ Assume for a moment that we set $D_{\sigma}=\{\epsi_2,~\epsi_3\}$ to be the support of $\sigma$. Then, of course, we have to say that $D_{\tau}=\{\epsi_1,~\epsi_4\}$ is the support of the involution $$\tau=\begin{pmatrix}1&2&3&4\\-1&2&3&-4\end{pmatrix}=s_{\epsi_1}s_{\epsi_4}=s_{\epsi_1-\epsi_4}s_{\epsi_1+\epsi_4},$$ and Theorem~\ref{theo:if_contained_then_Bruhat_B_n} fails immediately. Indeed, it follows from \cite[Proposition 2.1]{Ignatev1} that $\Omega_{\sigma}=\Omega_{s_{\epsi_2}}$ and $\Omega_{\tau}=\Omega_{s_{\epsi_1}}$, and once can easily check that $\Omega_{s_{\epsi_2}}\subseteq\overline{\Omega}_{s_{\epsi_1}}$ (and so $\Omega_{\sigma}\subseteq\overline{\Omega}_{\tau}$): if $$g=x_{\epsi_1-\epsi_2}(-t^{-1})h_{\epsi_1-\epsi_2}(t^{-1}),$$ then $g.f_{\epsi_1}\to f_{\epsi_2}$ as $t\to0$ (it is well-known that the Zariski closure of a constructive set coincides with its closure in the complex topology). At the same time, Theorem~\ref{mtheo:A_n_C_n} claims that $\sigma\nleq_B\tau$.

On the other hand, we may set $\{\epsi_2-\epsi_3,~\epsi_2+\epsi_3\}$ to be the support of $\sigma$. (According to \cite{Springer1}, this choice is ``more canonical'' than the previous one.) If we define the support of an involution in such a way, we neither have counterexamples to Theorem~\ref{theo:if_contained_then_Bruhat_B_n} no can prove it in general.}

\sst In\label{sst:if_contained_then_Bruhat_D_n} this subsection, we check that if $\sigma,~\tau\in\Bu(D_n)$, then $\Omega_{\sigma}$ is contained in the closure of $\Omega_{\tau}$. First, applying Lemma~\ref{lemm:Omega_union_Theta}, one can repeat the proof of Theorem~\ref{theo:if_contained_then_Bruhat_B_n} literally to show that $(R_{\sigma}^*)_{i,j}\leq(R_{\tau}^*)_{i,j}$ for all $i,~j$. According to Theorem~\ref{mtheo:A_n_C_n}, this means that $(R_{\sigma})_{i,j}\leq(R_{\tau})_{i,j}$ for all $i,~j$. Hence it remains to prove that the basis involutions $\sigma$ and $\tau$ satisfy the second condition in (\ref{formula:Bruhat_Dn}). To do this, we need to introduce some more notation.

Let $\Phi=D_n$ and $x\in\gt$. Given $r\leq n$ and two ordered $r$-tuples $P=\{p_1,\ldots,p_r\}$ and\break $Q=\{q_1,\ldots,q_r\}$, where $1\leq q_i\leq n$ and $p_i\in[\pm n]$ for $1\leq i\leq r$, we denote by $\Delta_P^Q(x)=\Delta_{p_1,\ldots,p_r}^{q_1,\ldots,q_r}(x)$ the minor of the matrix $x$ with the set of rows $P$ and the set of columns $Q$, i.e.,
\begin{equation*}
\Delta_P^Q(x)=\begin{vmatrix}
x_{p_1,q_1}&\ldots&x_{a_1,q_r}\\
\vdots&\ddots&\vdots\\
x_{p_r,q_1}&\ldots&x_{p_r,q_r}
\end{vmatrix}.
\end{equation*}
Given ordered tuples $I=\{i_1,\ldots,i_r\}$, $J=\{j_1,\ldots,j_s\}$, we denote by $I\cup J$ its concatenation, i.e., $I\cup J=\{i_1,\ldots,i_r,j_1,\ldots,j_s\}$.

If $P=\{p_1,\ldots,p_r\}$ is an $r$-tuple and $1\leq i<j\leq r$, then define the $r$-tuples $P^+[p_i,p_j]$ and $P^-[p_i,p_j]$ by the formula
\begin{equation*}
\begin{split}
P^+[p_i,p_j]&=\{p_1,\ldots,p_{i-1},p_{i+1},\ldots,p_{j-1},p_{j+1},\ldots,p_r,p_i,-p_i\},\\
P^-[p_i,p_j]&=\{p_1,\ldots,p_{i-1},p_{i+1},\ldots,p_{j-1},p_{j+1},\ldots,p_r,p_j,-p_j\}.\\
\end{split}
\end{equation*}

Next, if $P=\{p_1,\ldots,p_r\}$ is an $r$-tuple, $p_i\in[\pm n]$ for $1\leq i\leq r$, and $P'=\{p_{i_1},\ldots,p_{i_{2s}}\}$ is a tuple of even cardinality such that all $p_{i_j}$ are distinct and belong to $P$, we define the set $\So_{P,P'}$ of $r$-tuples by the following rule. If $P'=\varnothing$, then $\So_{P,P'}=\{P\}$. For $s\geq1$, we put $$\So_{P,P'}=\bigcup_{P_0\in\So_{P,P''}}\{P_0^+[p_{i_{2s-1}},p_{i_{2s}}],P_0^-[p_{i_{2s-1}},p_{i_{2s}}]\},$$ where, by definition, $P''=\{p_{i_1},\ldots,p_{i_{2s-2}}\}$. For example, if $P=\{1,2,-3,-4\}$, $P'=\{1,-3,2,-4\}$, then $\So_{P,P'}$ consists of four tuples: $\{1,-1,2,-2\},~\{1,-1,-4,4\},~\{-3,3,2,-2\}$, and $\{-3,3,-4,4\}$. Finally, if $P$, $Q$ are $r$-tuples and $P'$ is as above, we define the polynomial $$D_{P,P'}^Q(x)=\sum_{P_0\in\So_{P,P'}}\Delta_{P_0}^Q(x).$$

The following technical (but important) proposition is the key step in the proof of the main result of this subsection. Assume that $w\in\Bu(D_n)$ and $[a,-a]\times[b,-b]$, $a\geq b$, is an empty rectangle for~$w$. Let $P=\{p_1,\ldots,p_r\}$, $Q=\{q_1,\ldots,q_s\}$ be $r$-tuples, $1\leq q_i\leq b-1$, $p_i\in[\pm n]$ for all $i$, where $r=(R_w^*)_{a,b-1}$. Assume that $P=I\cup J\cup K$, where each element of $I$ (respectively, of $J$ and $K$) is from $a$ to $n$ (respectively, from $-n$ to $-a$ and from $-a+1$ to $-1$). Suppose that $|I\cup J|<n-a+1$ or
\begin{equation*}
\begin{split}
&\#\{\epsi_i-\epsi_j\in\Supp{w}\mid a\leq j\leq n\}\not\equiv|I|\pmod2,\\
&\#\{\epsi_i+\epsi_j\in\Supp{w}\mid a\leq j\leq n\}\not\equiv|J|\pmod2.\\
\end{split}
\end{equation*}
(For $|I\cup J|=n-a+1$, these two conditions are in fact equivalent.)

\propp{Let $P'=\{p_{i_1},p_{i_2},\ldots\}$ be a tuple\label{prop:zero_on_orbit} of even cardinality with distinct elements containing in $I\cup J$. Then $D_{P,P'}^Q(\lambda)=0$ for all $\lambda\in\Omega_w$.}{Pick $\lambda\in\Omega_w$. According to Lemma~\ref{lemm:Omega_union_Theta}, $\lambda=u.f_{w,\xi}$ for certain $u\in U$, $\xi\colon D\to\Cp^{\times}$, where $D=\Supp{w}$. Since $U$ is generated by $x_{\alpha}(s)$, $\alpha\in\Phi^+$, $s\in\Cp^{\times}$, there exist $\alpha_1,\ldots,\alpha_k\in\Phi^+$, $s_1,\ldots,s_k\in\Cp^{\times}$ such that $u=x_{\alpha_1(s_1)}\ldots x_{\alpha_k}(s_k)$. The proof is by induction on $k$. The base $k=0$ is trivial. Thus, we must prove the following fact: if $f\in\nt^*$, $\alpha\in\Phi^+$, $s\in\Cp^{\times}$, and all possible $D_{P_0,P_0'}^{Q_0}(f)$ are zero, then $D_{P,P'}^Q(x_{\alpha}(s).f)=0$.

Given $i,j\in[\pm n]$, $s\in\Cp$, consider the elementary transformations $r_{i,j}^s$ and $c_{i,j}^s$, where, for $f\in\nt^*$,
\begin{equation*}
\begin{split}
&(r_{i,j}^s(f))_{p,q}=\begin{cases}f_{i,q}+sf_{i,j}&\text{if }p=i>q,\\
f_{p,q}&\text{otherwise},
\end{cases}\\
&(c_{i,j}^s(f))_{p,q}=\begin{cases}f_{p,j}+sf_{q,j}&\text{if }p>q=j,\\
f_{p,q}&\text{otherwise}.
\end{cases}
\end{split}
\end{equation*}
Note that
\begin{equation*}
\begin{split}
&x_{\epsi_i-\epsi_j}(s).f=r_{i,j}^s(r_{-j,-i}^{-s}(c_{j,i}^s(c_{-i,-j}^{-s}(f)))),\\
&x_{\epsi_i+\epsi_j}(s).f=r_{i,-j}^s(r_{j,-i}^{-s}(c_{-i,j}^s(c_{-j,i}^{-s}(f)))).\\
\end{split}
\end{equation*}
Hence it is enough to prove that if all possible $D_{P_0,P_0'}^{Q_0}(f)$ are zero then $D_{P,P'}^Q(f')$ is zero, where $f'=r_{i,j}^s(r_{-j,-i}^{-s}(f))$, $r_{i,-j}^s(r_{j,-i}^{-s}(f))$, $c_{j,i}^s(c_{-i,-j}^{-s}(f))$ or $c_{-i,j}^s(c_{-j,i}^{-s}(f))$ for certain $i<j$, $s\in\Cp^{\times}$. We will consider this cases subsequently.

First, let $f'=r_{i,j}^s(r_{-j,-i}^{-s}(f))$. Given a tuple $T=\{t_1,t_2,\ldots\}$ with distinct elements and numbers $t_i\in T$, $t$, we denote by $T[t\to t_i]$ the tuple $T[t\to t_i]=\{t_1,\ldots,t_{i-1},t,t_{i+1},\ldots\}$. If $\pm i$, $\pm j\notin P$ (or $-i\in K$, $\pm j\notin P$, or $j\in I$, $j\notin P'$, $\pm i\notin P$, or $-i\in J$, $-i\notin P'$, $\pm j\notin P$, or $-i\in K$, $j\in I$, $j\notin P'$), then, obviously, $D_{P,P'}^Q(f')=D_{P,P'}^Q(f)=0$. If $i\in I$, $i\notin P'$, $\pm j\notin P$, then $|I[j\to i]\cup J|=|I\cup J|<n-a+1$, hence $$D_{P,P'}^Q(f')=D_{P,P'}^Q(f)+sD_{P[j\to i],P'}^Q(f)=0.$$ If $\pm i\notin P$, $-j\in J$, $-j\notin P'$ (respectively, $\pm i\notin P$, $-j\in K$), then $$D_{P,P'}^Q(f')=D_{P,P'}^Q(f)-sD_{P[-i\to -j],P'}^Q(f)=0,$$ because $|I\cup J[-i\to-j]|<n-a+1$ (respectively, because $P[-i\to\-j]=I\cup J\cup K[-i\to\-j]$).

Suppose $i\in I$, $i\in P'$, $\pm j\notin P$ (the case $-i\in J$, $-i\in P'$, $\pm j\notin P$ is similar). Given a tuple\break $T=\{t_1,t_2,\ldots\}$ with distinct elements and a number $t_i\in T$, we denote by $T\setminus\{t_i\}$ the tuple\break $T\setminus\{t_i\}=\{t_1,\ldots,t_{i-1},t_{i+1},\ldots\}$. If $t_i,~t_j,~\ldots$ are distinct elements of $T$, we define the tuple $T\setminus\{t_i,t_j,\ldots\}$ similarly. If $i=p_{i_{2k-1}}\in P'$ (respectively, $i=p_{i_{2k}}\in P'$) for some $k$, then denote $l=p_{i_{2k}}$ (respectively, $l=p_{i_{2k-1}}$). One has
$$D_{P,P'}^Q(f')=D_{P,P'}^Q(f)\pm sD_{(P[j\to l])[-i\to i],P'\setminus\{i,l\}}^Q(f)=0,$$
because $|((I\setminus\{i\})\cup(J\cup\{-i\}))[j\to l]|=|I\cup J|<n-a+1$.

Suppose now that $-j\in J$, $-j\in P'$, $\pm i\notin P$. (The cases $j\in I$, $j\in P'$, $\pm i\notin P$ and $-i\in J$, $-i\in P'$, $\pm j\notin P$ are similar.) If $-j=p_{i_{2k-1}}\in P'$ (respectively, $-j=p_{i_{2k}}\in P'$) for some $k$, then denote $l=p_{i_{2k}}$ (respectively, $l=p_{i_{2k-1}}$). Then
$$D_{P,P'}^Q(f')=D_{P,P'}^Q(f)\pm sD_{(P[i\to l])[j\to -j],P'\setminus\{i,l\}}^Q(f).$$
If $i\geq a$, then we can argue as in the previous paragraph. If $i<a$, then the last summand is zero because $|(I\cup(J\setminus\{-j\}))\setminus\{l\}|<|I\cup J|\leq n-a+1$.

If $i,~j\in I$, $i,~j\notin P'$, then one has
$$D_{P,P'}^Q(f')=D_{P,P'}^Q(f)+s\sum_{P_0\in\So_{P,P'}}\Delta_{P_0[j\to i]}^Q(f)=0,$$
because each minor in the last summand contains two same rows. The cases $-i,~-j\in J$, $-i,~-j\notin P'$, and $-i\in J$, $-i\in P'$, $j\in I$, $j\notin P'$, and $-i\in J$, $-i\notin P'$, $j\in I$, $j\in P'$, and $i,~j\in I$, $i\in P'$, $j\notin P'$, and $-i,~-j\in P'$, $-i\notin P'$, $-j\in P'$, and $-i\in K$, $j\in I$, $j\in P'$, and $-i\in K$, $-j\in J$, $-j\in P'$, and $-i\in K$, $-j\in J$, $-j\notin P'$, and $-i,~-j\in K$ are similar to this case.

Assume that $i,~j\in I$, $i\notin P'$, $j\in P'$. (The cases $i\in I$, $i\in P'$, $-j\in J$, $-j\notin P'$, and $i\in I$, $-j\in J$, $i\notin P'$, $-j\in P'$, and $-i,~-j\in J$, $-i\in P'$, $\-j\notin P'$ are similar.) We see that
\begin{equation*}
\begin{split}
D_{P,P'}^Q(f')&=D_{P,P'}^Q(f)+s\sum_{P_0\in\So_{P,P'},j\in P_0}\Delta_{P_0[j\to i]}^Q(f)
+s\sum_{P_0\in\So_{P,P'},j\notin P_0}\Delta_{P_0[j\to i]}^Q(f)\\
&-s\sum_{P_0\in\So_{P,P'},j\in P_0}\Delta_{P_0[-i\to-j]}^Q(f)-s^2\sum_{P_0\in\So_{P,P'},j\in P_0}\Delta_{(P_0[j\to i])[-i\to-j]}^Q(f).
\end{split}
\end{equation*}
Each minor in the second and the last summands contains two same rows. Thus,
$$D_{P,P'}^Q(f')=s\sum_{P_0\in\So_{P,P'},j\notin P_0}\Delta_{P_0[j\to i]}^Q(f)
-s\sum_{P_0\in\So_{P,P'},j\in P_0}\Delta_{P_0[-i\to-j]}^Q(f)=\pm sD_{P,P'[j\to i]}^Q(f)=0.$$

Assume now that $i,~j\in I$, $i~,j\in P'$. (The cases $-i~,-j\in J$, $-i,~-j\in P'$ and $i\in I$, $-j\in J$, $i,~-j\in P'$, and $-i\in J$, $j\in I$, $-i,~j\in P'$ are similar.) If $i=p_{i_{2k-1}}$ and $j=p_{i_{2k}}$, or $i=p_{i_{2k}}$ and  $j=p_{i_{2k-1}}$ for certain $k$, then
\begin{equation*}
\begin{split}
D_{P,P'}^Q(f')&=D_{P,P'}^Q(f)+s\sum_{P_0\in\So_{P,P'},i\in P_0}\Delta_{P_0[j\to i]}^Q(f)
-s\sum_{P_0\in\So_{P,P'},j\in P_0}\Delta_{P_0[-i\to-j]}^Q(f)\\
&=\pm s D_{P[-i\to i],P'\setminus\{i,j\}}^Q(f)\mp s D_{P[-i\to i],P'\setminus\{i,j\}}^Q(f)=0.
\end{split}
\end{equation*}
On the other hand, suppose that such a number $k$ does not exists. If $i=p_{i_{2k-1}}$ (respectively, $i=p_{i_{2k}}$) for certain $k$, then denote $i'=p_{i_{2k}}$ (respectively, $i'=p_{i_{2k-1}}$); define the number $j'$ in a similar way using $j$ instead of  $i$. Then
\begin{equation*}
\begin{split}
D_{P,P'}^Q(f')&=D_{P,P'}^Q(f)+s\sum_{P_0\in\So_{P,P'},i,j\in P_0}\Delta_{P_0[j\to i]}^Q(f)
-s\sum_{P_0\in\So_{P,P'},i,j\in P_0}\Delta_{P_0[-i\to-j]}^Q(f)\\
&+s\sum_{P_0\in\So_{P,P'},i,j'\in P_0}\Delta_{P_0[j\to i]}^Q(f)-
s\sum_{P_0\in\So_{P,P'},i',j\in P_0}\Delta_{P_0[-i\to-j]}^Q(f)\\
&-s^2\sum_{P_0\in\So_{P,P'},i,j\in P}\Delta_{P_0([j\to i])(-i\to-j)}^Q(f).
\end{split}
\end{equation*}
The second, the third and the last summands are zero, because each minor in these summands contains two same rows. Thus,
$$D_{P,P'}^Q(f')=s\sum_{P_0\in\So_{P,P'},i,j'\in P_0}\Delta_{P_0[j\to i]}^Q(f)
-s\sum_{P_0\in\So_{P,P'},i',j\in P_0}\Delta_{P_0[-i\to-j]}^Q(f)=\pm sD_{\wt P,\wt P'}^Q(f)=0,$$
where $\wt P=(P[-i\to i])[-j'\to j']$, $\wt P'=P'\setminus\{i,j\}$, because $\wt P$ satisfies all required conditions.

One of the most interesting cases is $i\in I$, $-j\in J$, $i,~-j\notin P'$. Here

\begin{equation*}
\begin{split}
D_{P,P'}^Q(f')&=D_{P,P'}^Q(f)+s\sum_{P_0\in\So_{P,P'}}\Delta_{P_0[j\to i]}^Q(f)
-s\sum_{P_0\in\So_{P,P'}}\Delta_{P_0[-i\to-j]}^Q(f)\\
&-s^2\sum_{P_0\in\So_{P,P'}}\Delta_{(P_0[j\to i])[-i\to-j]}^Q(f).
\end{split}
\end{equation*}
The last summand equals $\pm s^2D_{(P[j\to i])[-i\to-j],P'}^Q(f)=0$. Hence
\begin{equation*}
D_{P,P'}^Q(f')=s\sum_{P_0\in\So_{P,P'}}\Delta_{P_0[j\to i]}^Q(f)
-s\sum_{P_0\in\So_{P,P'}}\Delta_{P_0[-i\to-j]}^Q(f)
=\pm s D_{P,P'\cup\{i,-j\}}^Q(f)=0.
\end{equation*}
The case $-i\in J$, $j\in I$, $-i,~j\notin P'$ is completely similar.

The cases $i\in I$, $i\notin P'$, $-j\in K$, and $i\in I$, $i\notin P'$, $-j\in K$, and $-i\in J$, $-i\notin P'$, $-j\in K$, and $-i\in I$, $i\in P'$, $-j\in K$ can not occur because of the definition of $I$, $J$, $K$. Thus, we have considered all possible cases for $f'=r_{i,j}^s(r_{-j,-i}^{-s}(f))$.

Second, let $f'=r_{i,-j}^s(r_{j,-i}^{-s}(f))$. The proof in this case is similar to the proof for $f'=r_{i,j}^s(r_{-j,-i}^{-s}(f))$, so we skip the details. Next, let $f'=c_{j,i}^s(c_{-i,-j}^{-s}(f))$. If $j\notin Q$ then, clearly, $$D_{P,P'}^Q(f')=D_{P,P'}^Q(f)=0.$$ If $j\in Q$, $i\notin Q$, then $$D_{P,P'}^Q(f')=D_{P,P'}^Q(f)+sD_{P,P'}^{Q[i\to j]}(f)=0.$$ If $i,~j\in Q$, then
\begin{equation*}
D_{P,P'}^Q(f')=D_{P,P'}^Q(f)+s\sum_{P_0\in\So_{P,P'}}\Delta_{P_0}^{Q[i\to j]}(f)=0,
\end{equation*}
because each minor in the last summand contains two same columns. Finally, let $f'=c_{-i,j}^s(c_{-j,i}^{-s}(f))$, then $D_{P,P'}^Q(f')=D_{P,P'}^Q(f)=0$. The proof is complete.}

Now we can prove the main result of this subsection.\newpage

\theop{Let $\sigma,~\tau\in\Bu(D_n)$. If\label{theo:if_contained_then_Bruhat_D_n} $\Omega_{\sigma}\subseteq\overline{\Omega}_{\tau}$\textup, then $\sigma\leq_B\tau$.}{Assume that $\sigma\nleq_B\tau$. If there exist $i,~j$ such that $(R_{\sigma}^*)_{i,j}>(R_{\tau}^*)_{i,j}$, then, using Lemma~\ref{lemm:Omega_union_Theta} and arguing as in the proof of Theorem~\ref{theo:if_contained_then_Bruhat_B_n}, one can check that $\Omega_{\sigma}\not\subseteq\overline{\Omega}_{\tau}$. Hence $R_{\sigma}^*\leq R_{\tau}^*$, and, by Theorem~\ref{mtheo:A_n_C_n}, $R_{\sigma}\leq R_{\tau}$. Thus, according to (\ref{formula:Bruhat_Dn}), there exist $a$, $b$ such that $[a,-a]\times[b,-b]$ is an empty rectangle for $\sigma$ and $\tau$, $(R_{\sigma})_{-(a-1),b-1}=(R_{\tau})_{-(a-1),b-1}$, but $(R_{\sigma})_{-(a-1),n}\not\equiv(R_{\tau})_{-(a-1),n}\pmod2$. Since $\sigma$ and $\tau$ are involutions, we may assume without loss of generality that $a\geq b$.

Recall the notion of $X_w$ for $w\in W$ from Subsection\ref{sst:B_n_D_n_definitions}. Given an arbitrary element $w\in W$, $-n\leq p\leq q\leq-1$ and $1\leq r\leq s\leq n$, denote
\begin{equation*}
w_{p,q}^{r,s}=\#\{(i,j)\mid p\leq i\leq q,~r\leq j\leq s\text{ and }(X_w)_{i,j}=1\}.
\end{equation*}
In other words, $w_{p,q}^{r,s}$ is the number of rooks in $X_w$ which rows (respectively, columns) are between $p$ and $q$ (respectively, between $r$ and $s$). By definition of $R_w$,
\begin{equation*}
(R_w)_{-(a-1),-n}=w_{-(b-1),-1}^{1,b-1}+w_{-(a-1),-b}^{1,b-1}+w_{-(b-1),-1}^{b,a-1}+w_{-(a-1),-b}^{b,a-1}+w_{-(b-1),-1}^{a,n}+w_{-(a-1),-b}^{a,n}.
\end{equation*}
If $w$ is an involution in the Weyl group $W$, then $w_{-n,-a}^{1,b-1}=w_{-(b-1),-1}^{a,n}$, while the numbers $w_{-(b-1),-1}^{1,b-1}$ and $w_{-(a-1),-b}^{b,a-1}$ are even. Furthermore, in this case,
\begin{equation*}
w_{-(a-1),-1}^{1,b-1}=w_{-(b-1),-1}^{1,b-1}+w_{-(a-1),-b}^{1,b-1}=w_{-(b-1),-1}^{1,a-1}=w_{-(b-1),-1}^{1,b-1}+w_{-(b-1),-1}^{b,a-1},
\end{equation*}
so $w_{-(a-1),-b}^{1,b-1}=w_{-(b-1),-1}^{b,a-1}$. If, in addition, $[-a,a]\times[-b,b]$ is an empty rectangle for $w$, then $w_{-(a-1),-b}^{a,n}=0$, thus, $(R_w)_{-(a-1),b-1}\equiv w_{-n,-a}^{1,b-1}\pmod2$.

By our assumptions, $(R_{\sigma})_{-(a-1),n}\not\equiv(R_{\tau})_{-(a-1),n}\pmod2$, hence $\sigma_{-n,-a}^{1,b-1}\not\equiv\tau_{-n,-a}^{1,b-1}\pmod2$ (this is the key observation in the proof). Denote $P=I\cup J\cup K$, $Q=\{1.\ldots,b-1\}$, where
\begin{equation*}
\begin{split}
I&=\{i\in\{a,\ldots,n\}\mid\sigma(i)\in Q\},\\
J&=\{j\in\{-n,\ldots,-a\}\mid\sigma(j)\in Q\},\\
K&=\{k\in\{-(a-1),\ldots,-1\}\mid\sigma(k)\in Q\}.
\end{split}
\end{equation*}
Note that $|K|=(R_{\sigma})_{-(a-1),b-1}=(R_{\tau})_{-(a-1),b-1}$, while $|I|+|J|=n-a+1$ (the latter equality follows from the fact that $[-a,a]\times[-b,b]$ is an empty rectangle for $\sigma$). Moreover, it follows from $\sigma_{-n,-a}^{1,b-1}\not\equiv\tau_{-n,-a}^{1,b-1}\pmod2$ that $\#\{\epsi_i-\epsi_j\in\Supp{\tau}\mid a\leq j\leq n\}\not\equiv|I|\pmod2$ (or, equivalently, $\#\{\epsi_i+\epsi_j\in\Supp{w}\mid a\leq j\leq n\}\not\equiv|J|\pmod2$). Hence, accordingly to Proposition~\ref{prop:zero_on_orbit}, $D_{P,P'}^Q(\lambda)=0$ for all $\lambda\in\Omega_{\tau}$ and all subsets $P'\subset P$ of even cardinality. In particular, $$D_{P,\varnothing}^Q(\lambda)=\Delta_P^Q(\lambda)=0$$ for all $\lambda\in\Omega_{\tau}$ (and, consequently, for all $\lambda\in\overline{\Omega}_{\tau}.$) But $\Delta_P^Q(f_{\sigma})=\pm1$ by definition of $f_{\sigma}$. Thus, $\Omega_{\sigma}\notin\overline{\Omega}_{\tau}$. This concludes the proof.}

\sst \label{sst:if_Bruhat_then_contained} In this subsection we present a conjectural way how to prove that if $\sigma,~\tau\in\Bu(\Phi)$ and $\sigma\leq_B\tau$ then $\Omega_{\sigma}\subseteq\overline{\Omega}_{\tau}$. (For simplicity, we consider only the case $\Phi=B_n$ and give some additional remarks for $D_n$ at the end of the subsection.) To do this, we need to describe the covering relation on $\Bu(\Phi)$ with respect to the Bruhat order. The covering relation on $\Iu(\Phi)$ was described by F. Incitti in \cite{Incitti3}. We will state a corollary of his description for $C_n$ in appropriate terms. To each involution $w\in\Iu(C_n)$ we assign the number $$d(\sigma)=|\{i\in\{1,~\ldots,~n\}\mid2\epsi_i\in\Supp{\sigma}\}|.$$ Let $\sigma,~\tau\in\Iu(C_n)$, and $D_{\sigma}=\Supp{\sigma}$, $D_{\tau}=\Supp{\tau}$. (In fact, we will apply this corollary to $B_n$, but it is more convenient to formulate it for $C_n$, because we will use the notion of the support of an arbitrary involution.) In the tables below we consider certain special cases of ``relative positions'' of $D_{\sigma}$ and $D_{\tau}$ in $\Phi^+$, which are needed to formulate this corollary.

\newpage\begin{center}
\textbf{Table 1. Case $d(\sigma)=d(\tau)$, first part}
\end{center}
\vspace{-0.7cm}\begin{longtable}{||l|l|l||l|l|l||l|l|l||}
\hline
&$D_{\sigma}\setminus D_{\tau}$&$D_{\tau}\setminus D_{\sigma}$&&$D_{\sigma}\setminus D_{\tau}$&$D_{\tau}\setminus D_{\sigma}$
&&$D_{\sigma}\setminus D_{\tau}$&$D_{\tau}\setminus D_{\sigma}$\\
\hline\hline1
1
&\begin{tabular}[t]{l}$\epsi_i+\epsi_j$,\\
$i<k<j$
\end{tabular}&
\begin{tabular}[t]{l}$\epsi_i+\epsi_k$
\end{tabular}&
2
&\begin{tabular}[t]{l}$\epsi_k-\epsi_j$,\\
$\epsi_i+\epsi_l$,\\
$i<k<j<l$
\end{tabular}&
\begin{tabular}[t]{l}$\epsi_i+\epsi_j$,\\
$\epsi_k-\epsi_l$
\end{tabular}&
3
&\begin{tabular}[t]{l}$\epsi_i-\epsi_j$,\\
$\epsi_k+\epsi_l$,\\
$i<k<j<l$
\end{tabular}&
\begin{tabular}[t]{l}$\epsi_i-\epsi_l$,\\
$\epsi_k+\epsi_j$
\end{tabular}\\
\hline
4
&\begin{tabular}[t]{l}$\epsi_i-\epsi_j$,\\
$i<j$
\end{tabular}&
\begin{tabular}[t]{l}$\epsi_i+\epsi_j$
\end{tabular}&
5
&\begin{tabular}[t]{l}$\epsi_i+\epsi_j$,\\
$\epsi_k+\epsi_l$,\\
$i<k<j<l$
\end{tabular}&
\begin{tabular}[t]{l}$\epsi_i+\epsi_k$,\\
$\epsi_j+\epsi_l$
\end{tabular}&
6
&\begin{tabular}[t]{l}$\epsi_i+\epsi_j$,\\
$\epsi_k-\epsi_l$,\\
$i<k<j<l$
\end{tabular}&
\begin{tabular}[t]{l}$\epsi_i+\epsi_k$,\\
$\epsi_j-\epsi_l$
\end{tabular}\\
\hline
7
&\begin{tabular}[t]{l}$\epsi_i+\epsi_l$,\\
$\epsi_k-\epsi_j$,\\
$i<k<j<l$
\end{tabular}&
\begin{tabular}[t]{l}$\epsi_i+\epsi_k$
\end{tabular}&
8
&\begin{tabular}[t]{l}$\epsi_k-\epsi_j$,\\
$i<k<j$
\end{tabular}&
\begin{tabular}[t]{l}$\epsi_i-\epsi_j$
\end{tabular}&
9
&\begin{tabular}[t]{l}$\epsi_k+\epsi_j$,\\
$i<k<j$
\end{tabular}&
\begin{tabular}[t]{l}$\epsi_i+\epsi_j$
\end{tabular}\\
\hline
10
&\begin{tabular}[t]{l}$\epsi_i-\epsi_k$,\\
$i<k<j$
\end{tabular}&
\begin{tabular}[t]{l}$\epsi_i-\epsi_j$
\end{tabular}&
11
&\begin{tabular}[t]{l}$\epsi_i-\epsi_k$,\\
$\epsi_j-\epsi_l$,\\
$i<k<j<l$
\end{tabular}&
\begin{tabular}[t]{l}$\epsi_i-\epsi_j$,\\
$\epsi_k-\epsi_l$
\end{tabular}&
12
&\begin{tabular}[t]{l}$\epsi_i-\epsi_k$,\\
$\epsi_j+\epsi_l$,\\
$i<k<j<l$
\end{tabular}&
\begin{tabular}[t]{l}$\epsi_i-\epsi_j$,\\
$\epsi_k+\epsi_l$
\end{tabular}\\
\hline
13
&\begin{tabular}[t]{l}$\epsi_i+\epsi_l$,\\
$\epsi_k+\epsi_j$,\\
$i<k<j<l$
\end{tabular}&
\begin{tabular}[t]{l}$\epsi_i+\epsi_j$,\\
$\epsi_k+\epsi_l$
\end{tabular}&
14
&\begin{tabular}[t]{l}$\epsi_i-\epsi_j$,\\
$\epsi_k-\epsi_l$,\\
$i<k<j<l$
\end{tabular}&
\begin{tabular}[t]{l}$\epsi_i-\epsi_l$,\\
$\epsi_k-\epsi_j$
\end{tabular}&
15
&\begin{tabular}[t]{l}$\epsi_i-\epsi_l$,\\
$\epsi_k+\epsi_j$,\\
$i<k<j<l$
\end{tabular}&
\begin{tabular}[t]{l}$\epsi_i+\epsi_j$,\\
$\epsi_k-\epsi_l$
\end{tabular}\\
\hline
16
&\begin{tabular}[t]{l}$\epsi_i-\epsi_j$,\\
$\epsi_k+\epsi_l$,\\
$i<k<j<l$
\end{tabular}&
\begin{tabular}[t]{l}$\epsi_i+\epsi_l$,\\
$\epsi_k-\epsi_j$
\end{tabular}&
17
&\begin{tabular}[t]{l}$\epsi_i-\epsi_k$,\\
$\epsi_j-\epsi_l$,\\
$i<k<j<l$
\end{tabular}&
\begin{tabular}[t]{l}$\epsi_i-\epsi_l$
\end{tabular}&
18
&\begin{tabular}[t]{l}$\epsi_i-\epsi_k$,\\
$\epsi_j+\epsi_l$,\\
$i<k<j<l$
\end{tabular}&
\begin{tabular}[t]{l}$\epsi_i+\epsi_l$
\end{tabular}\\
\hline
19
&\begin{tabular}[t]{l}$\varnothing$,\\
$i<j$
\end{tabular}&
\begin{tabular}[t]{l}$\epsi_i-\epsi_j$
\end{tabular}&&&&&&\\
\hline
\end{longtable}

\begin{center}
\textbf{Table 2. Case $d(\sigma)=d(\tau)$, second part}
\end{center}
\vspace{-0.7cm}\begin{longtable}{||l|l|l||l|l|l||l|l|l||}
\hline
&$D_{\sigma}\setminus D_{\tau}$&$D_{\tau}\setminus D_{\sigma}$&&$D_{\sigma}\setminus D_{\tau}$&$D_{\tau}\setminus D_{\sigma}$
&&$D_{\sigma}\setminus D_{\tau}$&$D_{\tau}\setminus D_{\sigma}$\\
\hline\hline
1
&\begin{tabular}[t]{l}$2\epsi_j$,\\
$i<j$
\end{tabular}&
\begin{tabular}[t]{l}$2\epsi_i$
\end{tabular}&
2
&\begin{tabular}[t]{l}$\epsi_i+\epsi_j$, $2\epsi_k$,\\
$i<k<j$
\end{tabular}&
\begin{tabular}[t]{l}$2\epsi_i$, $\epsi_k+\epsi_j$\end{tabular}&
3
&\begin{tabular}[t]{l}$\epsi_i-\epsi_j$, $2\epsi_k$,\\
$i<k<j$
\end{tabular}&
\begin{tabular}[t]{l}$2\epsi_i$, $\epsi_k-\epsi_j$\end{tabular}\\
\hline
4
&\begin{tabular}[t]{l}$\epsi_i-\epsi_k$, $2\epsi_j$,\\
$i<k<j$
\end{tabular}&
\begin{tabular}[t]{l}$2\epsi_i$\end{tabular}&
5
&\begin{tabular}[t]{l}$\epsi_i-\epsi_k$, $2\epsi_j$,\\
$i<k<j$
\end{tabular}&
\begin{tabular}[t]{l}$\epsi_i-\epsi_j$, $2\epsi_k$\end{tabular}&
6
&\begin{tabular}[t]{l}$\epsi_i+\epsi_j$, $2\epsi_k$,\\
$i<k<j$
\end{tabular}&
\begin{tabular}[t]{l}$\epsi_i+\epsi_k$, $2\epsi_j$\end{tabular}\\
\hline
\end{longtable}

\begin{center}
\textbf{Table 3. Case $d(\sigma)<d(\tau)$}
\end{center}
\vspace{-0.7cm}\begin{longtable}{||l|l|l||l|l|l||l|l|l||}
\hline
&$D_{\sigma}\setminus D_{\tau}$&$D_{\tau}\setminus D_{\sigma}$&&$D_{\sigma}\setminus D_{\tau}$&$D_{\tau}\setminus D_{\sigma}$
&&$D_{\sigma}\setminus D_{\tau}$&$D_{\tau}\setminus D_{\sigma}$\\
\hline\hline
1
&\begin{tabular}[t]{l}$\emptyset$\end{tabular}&
\begin{tabular}[t]{l}$2\epsi_i$\end{tabular}&
2
&\begin{tabular}[t]{l}$\epsi_i+\epsi_j$,\\
$i<j$
\end{tabular}&
\begin{tabular}[t]{l}$2\epsi_i$, $2\epsi_j$\end{tabular}&
3
&\begin{tabular}[t]{l}$\epsi_i-\epsi_j$,\\
$i<j$
\end{tabular}&
\begin{tabular}[t]{l}$2\epsi_i$\end{tabular}\\
\hline
4
&\begin{tabular}[t]{l}$\epsi_i+\epsi_j$,\\
$\epsi_k-\epsi_l$
\end{tabular}&
\begin{tabular}[t]{l}$2\epsi_i$, $\epsi_k+\epsi_j$,\\
$i<k<j<l$\end{tabular}&
5
&\begin{tabular}[t]{l}$\epsi_i+\epsi_l$,\\
$\epsi_k-\epsi_j$
\end{tabular}&
\begin{tabular}[t]{l}$2\epsi_i$, $\epsi_k+\epsi_l$,\\
$i<k<j<l$\end{tabular}&
6
&\begin{tabular}[t]{l}$\epsi_i-\epsi_l$,\\
$\epsi_k-\epsi_j$
\end{tabular}&
\begin{tabular}[t]{l}$2\epsi_i$, $\epsi_k-\epsi_l$,\\
$i<k<j<l$\end{tabular}\\
\hline
\end{longtable}

\begin{center}
\textbf{Table 4. Case $d(\sigma)>d(\tau)$}
\end{center}
\vspace{-0.7cm}\begin{longtable}{||l|l|l||l|l|l||l|l|l||}
\hline
&$D_{\sigma}\setminus D_{\tau}$&$D_{\tau}\setminus D_{\sigma}$&&$D_{\sigma}\setminus D_{\tau}$&$D_{\tau}\setminus D_{\sigma}$
&&$D_{\sigma}\setminus D_{\tau}$&$D_{\tau}\setminus D_{\sigma}$\\
\hline\hline
1
&\begin{tabular}[t]{l}$\epsi_i-\epsi_j$, $2\epsi_k$,\\
$i<k<j$
\end{tabular}&
\begin{tabular}[t]{l}$\epsi_i+\epsi_k$\end{tabular}&
2
&\begin{tabular}[t]{l}$2\epsi_j$
\end{tabular}&
\begin{tabular}[t]{l}$\epsi_i+\epsi_j$,\\
$i<j$
\end{tabular}&
3
&\begin{tabular}[t]{l}$2\epsi_k,2\epsi_j$
\end{tabular}&
\begin{tabular}[t]{l}$\epsi_i+\epsi_k$,\\
$i<k<j$
\end{tabular}\\
\hline
4
&\begin{tabular}[t]{l}$\epsi_i-\epsi_k$, $2\epsi_j$,\\
$i<k<j$
\end{tabular}&
\begin{tabular}[t]{l}$\epsi_i+\epsi_j$
\end{tabular}&
5
&\begin{tabular}[t]{l}$\epsi_i-\epsi_k$,\\
$2\epsi_j$, $2\epsi_l$,\\
$i<k<j<l$
\end{tabular}&
\begin{tabular}[t]{l}$\epsi_i-\epsi_l$,\\
$\epsi_k+\epsi_j$
\end{tabular}&
6
&\begin{tabular}[t]{l}$\epsi_i+\epsi_l$,\\
$2\epsi_k$, $2\epsi_j$,\\
$i<k<j<l$
\end{tabular}&
\begin{tabular}[t]{l}$\epsi_i+\epsi_k$,\\
$\epsi_j+\epsi_l$
\end{tabular}\\
\hline
7
&\begin{tabular}[t]{l}$\epsi_i-\epsi_l$,\\
$2\epsi_k$, $2\epsi_j$,\\
$i<k<j<l$
\end{tabular}&
\begin{tabular}[t]{l}$\epsi_i+\epsi_k$,\\
$\epsi_j-\epsi_l$
\end{tabular}&
8
&\begin{tabular}[t]{l}$\epsi_i-\epsi_j$,\\
$2\epsi_k$, $2\epsi_l$,\\
$i<k<j<l$
\end{tabular}&
\begin{tabular}[t]{l}$\epsi_i+\epsi_k$
\end{tabular}&
9
&\begin{tabular}[t]{l}$\epsi_i-\epsi_k$,\\
$2\epsi_j$, $2\epsi_l$,\\
$i<k<j<l$
\end{tabular}&
\begin{tabular}[t]{l}$\epsi_i+\epsi_j$
\end{tabular}\\
\hline
\end{longtable}

Given involutions $\sigma$, $\tau\in\Iu(C_n)$, we say that the pair $(\tau,\sigma)$ is \emph{of type} (or, equivalently, \emph{belongs to case}) $a.b$ if $D_{\sigma}=\Supp{\sigma}$ and $D_{\tau}=\Supp{\tau}$ are as in the case $b$ in Table $a$ (for some $i,~k,~j,~l$). We also say that $(\tau,\sigma)$ is an \emph{admissible pair} if it is of type $a.b$ for certain $a,~b$. From \cite[Chapter 6]{Incitti1} (see also \cite[p. 76--91]{Incitti3}) one can immediately deduce the following

\mcoro{Let\label{coro:C_n_covering_relation} $\sigma,~\tau\in\Iu(C_n)$ and $\sigma<_B\tau$. Then there exist $\sigma_1,~\ldots,~\sigma_k\in\Iu(C_n)$ such that $\sigma_1=\tau$\textup, $\sigma_k=\sigma$ and $(\sigma_i,\sigma_{i+1})$ is an admissible pair for all $1\leq i\leq k-1$.}

Actually, this corollary does not describe the covering relation on $\Iu(C_n)$ (there are some additional conditions on $\sigma$ and $\tau$), but we will use only this part of the Incitti's description.

Since the Weyl groups of $B_n$ and $C_n$ coincide, we have the notion of a basis involution in $W(C_n)$; we denote the set of all basis involutions in $W(C_n)$ by $\Bu(C_n)$. We say that a pair $(\tau,\sigma)$ of basis involutions from $\Iu(C_n)$ is \emph{basis-admissible} if it is of type $1.b$ for certain $b$.

\hypo{Let\label{conj:basis_admissible} $\sigma,~\tau\in\Bu(C_n)$ and $\sigma<_B\tau$. Then there exist $\sigma_1,~\ldots,~\sigma_k\in\Bu(C_n)$ such that $\sigma_1=\tau$\textup, $\sigma_k=\sigma$ and $(\sigma_i,\sigma_{i+1})$ is basis-admissible for all $1\leq i\leq k-1$.}

We checked that this conjecture is true for $n\leq7$. It is easy to see that if this conjecture is true for all $n$ that $\sigma\leq_B\tau$ implies $\Omega_{\sigma}\subseteq\overline{\Omega}_{\tau}$. Indeed, we may assume without loss of generality that $(\tau,\sigma)$ is a basis-admissible pair. For basis-admissible pairs, the proof is case-by-case. For example, suppose that $(\tau,\sigma)$ is of type 1.12, i.e.,
\begin{equation*}
\begin{split}
&\Supp{\sigma}\setminus\Supp{\tau}=\{\epsi_i-\epsi_k,~\epsi_j+\epsi_l\},\\
&\Supp{\tau}\setminus\Supp{\sigma}=\{\epsi_i-\epsi_j,~\epsi_k+\epsi_j\}
\end{split}
\end{equation*}
for some $1\leq i<k<j<l\leq n$. Put $g(s)=x_{\epsi_k-\epsi_j}(t^{-2})h_{\epsi_i-\epsi_j}(t^{-1})h_{\epsi_k+\epsi_l}(-It^{-1})$
and $f=g(s).f_{\tau}$. (Here $I$ is the usual imaginary unit.) One can easily check by straightforward matrix calculations that, given $\alpha\in\Phi^+$,
\begin{equation*}
f(e_{\alpha})=\begin{cases}1,&\text{if either }\alpha=\epsi_i-\epsi_k\text{ or }\alpha=\epsi_j+\epsi_l,\\
t^2,&\text{if }\alpha=\epsi_i-\epsi_j,\\
-t^2,&\text{if }\alpha=\epsi_k+\epsi_l,\\
f_{\tau}(e_{\alpha})&\text{otherwise}.\\
\end{cases}
\end{equation*}
Thus, $f\to f_{\sigma}$ as $t\to0$. All other cases can be considered similarly. (For $\Phi=D_n$, one should use \cite[Chapter 7]{Incitti1} instead of \cite[Chapter 6]{Incitti1} for the description of the covering relation on $\Iu(D_n)$.)

\exam{Note that\label{exam:if_Bruhat_not_contained_B_n} $\sigma\leq_B\tau$ does not imply $\Omega_{\sigma}\subseteq\overline{\Omega}_{\tau}$ for non-basis involutions in $\Iu(B_n)$. Indeed, let $n=4$ and
\begin{equation*}
\sigma=\begin{pmatrix}1&2&3&4\\-4&2&3&-1\end{pmatrix}=s_{\epsi_1+\epsi_4},~
\tau=\begin{pmatrix}1&2&3&4\\-1&-3&-2&4\end{pmatrix}=s_{\epsi_1}s_{\epsi_2+\epsi_3}.
\end{equation*}}
Clearly, there is the only possibility to define the supports of these involutions, precisely,\break $\Supp{\sigma}=\{\epsi_1+\epsi_4\}$ and $\Supp{\tau}=\{\epsi_1,~\epsi_2+\epsi_3\}$. One can immediately deduce from Theorem~\ref{mtheo:A_n_C_n} that $\sigma<_B\tau$.

On the other hand, by definition of the coadjoint action, if $x,~y\in\nt$, $f\in\nt^*$, then
$$((\exp{x}).f)(y)=f(y)-f([x,y])+\dfrac{1}{2}f([x,[x,y]])-\ldots=f(\exp\ad{-x}(y)),$$
where, as usual, $\ad{x}(y)=[x,y]$. This implies that, given $\alpha\in\Phi^+$, if there are no $\beta\in\Phi^+$ such that $\alpha\leq\beta$ with respect to the natural order on $\Phi$ (i.e., $\beta-\alpha$ is zero or a sum of positive roots), then $\lambda(e_{\alpha})=0$ for all $\lambda\in\Theta_f$. Thus, $\lambda(e_{\epsi_1+\epsi_4})=0$ for all $\lambda\in\Omega_\tau$, but $f_{\sigma}(\epsi_1+\epsi_4)\neq0$, so $\Omega_{\sigma}\not\subset\overline{\Omega}_{\tau}$.

\newpage\sect{Concluding\label{sect:concluding_remarks} remarks}

\sst\label{sst:dim_Omega} Let $w$ be an involution from $\Bu(\Phi)$, where $\Phi=B_n$ or $D_n$. Being an orbit of a connected unipotent group on an affine variety $\nt^*$, $\Omega_w$ is a closed subvariety of $\nt^*$. In this subsection, we present a formula for the dimension of $\Omega_w$ (cf. \cite[Proposition 4.1]{Ignatev3}, \cite[Theorem 3.1]{Ignatyev4}). Recall the definition of the length $l(w)$ of an element $w\in W$.

\medskip\theop{Let $\Phi=B_n$ or $D_n$\textup, and $w\in\Bu(\Phi)$.\label{theo:dim_Omega} One has \begin{equation*}\dim\Omega_{w}=l(w).\vspace{-0.5cm}\end{equation*}}{Denote $D=\Supp{w}$. We claim that if $\xi_1$ and $\xi_2$ are two distinct maps from the set $D$ to $\Cp^{\times}$, then $$\Theta_{w,\xi_1}\neq\Theta_{w,\xi_2}.$$ Indeed, for $\Phi=B_n$ (respectively, $D_n$) let $\wt U$ be the group of all $(2n+1)\times(2n+1)$ (respectively, $2n\times2n$) upper-triangular matrices with $1$'s on the diagonal. Since $w$ is an involution in $S_{\pm n}$,\break \cite[Theorem 1.4]{Panov} implies that $\wt\Theta_{w,\xi_1}\neq\wt\Theta_{w,\xi_2}$, where $\wt\Theta_{w,\xi_1}$ and $\wt\Theta_{w,\xi_2}$ denote the respective $\wt U$-orbits of $f_{w,\xi_1}$ and $f_{w,\xi_2}$ under the coadjoint action of $\wt U$ on the space of all lower-triangular matrices with zeroes on the diagonal (see Subsection~\ref{sst:A_n_C_n_definitions} for the definitions). Since $U\subseteq\wt U$, one has $\Theta_{w,\xi_1}\subseteq\wt\Theta_{w,\xi_1}$ and $\Theta_{w,\xi_2}\subseteq\wt\Theta_{w,\xi_2}$, hence $\Theta_{w,\xi_1}\neq\Theta_{w,\xi_2}$, as required.

Let $Z_B=\Stab_Bf_{w}$ be the stabilizer of $f_{w}$ in $B$. One has $$\dim\Omega_{w}=\dim B-\dim Z_B.$$
Recall that $B\cong U\rtimes H$ as algebraic groups. It was shown in the proof of Lemma~\ref{lemm:Omega_union_Theta} that if $h\in H$, then there exists $\xi\colon D\to\Cp^{\times}$ such that $h.f_{w}=f_{w,\xi}$. Hence if $g=uh\in Z_B$, then $$f_{w}=(uh).f_{w}=u.f_{w,\xi},$$
so $f_{w}\in\Theta_{w,\xi}$. In follows from the first paragraph of the proof that $f_{w}=f_{w,\xi}$. This means that the map $$Z_U\times Z_H\to Z_B\colon (u,h)\mapsto uh$$ is an isomorphism of algebraic varieties, where $Z_U=\Stab_Uf_{w}$ (respectively, $Z_H=\Stab_Hf_{w}$) is the stabilizer of $f_{w}$ in $U$ (respectively, in $H$). Hence $$\dim Z_B=\dim Z_U+\dim Z_H.$$

Let $\Theta_w$ be the $U$-orbit of $f_w$. By \cite[Theorem 1.2]{Ignatev1}, $\dim\Theta_w=l(w)-|D|$, so $$\dim Z_U=\dim U-\dim\Theta_w=\dim U-l(w)+|D|.$$ On the other hand, put $X=\bigcup_{\xi\colon D\to\Cp^{\times}}\{f_{w,\xi}\}$. It follows from Lemma~\ref{lemm:Omega_union_Theta} and the first paragraph of the proof that $X=\{h.f_{w},h\in H\}$, the $H$-orbit of $f_{w}$. Consequently, $$\dim Z_H=\dim H-\dim X=\dim H-|D|,$$ because $X$ is isomorphic as affine variety to the product of $|D|$ copies of $\Cp^{\times}$. Thus,
\begin{equation*}
\begin{split}
\dim\Omega_{w}&=\dim B-\dim Z_B=(\dim U+\dim H)-(\dim Z_U+\dim Z_H)\\
&=\dim U+\dim H-(\dim U-l(w)-|D|)-(\dim H-|D|)=l(w).
\end{split}
\end{equation*} The proof is complete.}

\sst\label{sst:Schubert} In the remainder of the paper, we briefly discuss a conjectural geometrical approach to orbits associated with involutions in terms of tangent cones to Schubert varieties. Recall that $W$ is isomorphic to $N_G(H)/H$, where $N_G(H)$ is the normalizer of $H$ in $G$. The \emph{flag variety} $\Fo=G/B$ can be decomposed into the union $\Fo=\bigcup_{w\in W}\Xu_w^{\circ}$, where $\Xu_w^{\circ}=B\dot wB/B$
is called the \emph{Schubert cell} corresponding to $w$. (Here $\dot w$ is a representative of $w$ in $N_G(H)$.) By definition, the \emph{Schubert variety} $\Xu_w$ is the closure of
$\Xu_w^{\circ}$ in~$\Fo$ with respect to Zariski topology. Note that $p=\Xu_{\id}=B/B$ is contained in $\Xu_w$ for all $w\in W$. One has $\Xu_w\subseteq\Xu_{w'}$ if and only if $w\leq_Bw'$. Let $T_w$ be the tangent space and $C_w$ the tangent cone to~$\Xu_w$ at
the point $p$ (see \cite{BileyLakshmibai1} for detailed constructions); by definition, $C_w\subseteq T_w$,
and if $p$ is a~regular point of~$\Xu_w$, then $C_w=T_w$. Of course, if $w\leq_Bw'$, then $C_w\subseteq C_{w'}$.

Let $T=T_p\Fo$ be the tangent space to $\Fo$ at $p$. It can be naturally identified with $\nt^*$ by the following way: since $\Fo=G/B$, $T$ is isomorphic to the factor $\gt/\bt\cong\nt^*$. Next, $B$ acts on $\Fo$ by conjugation. Since $p$ is invariant under this action, the action on $T=\nt^*$ is induced. One can check that this action coincides with the action of $B$ on $\nt^*$ defined above. The tangent cone $C_w\subseteq T_w\subseteq T=\nt^*$
is $B$-invariant, so it splits into a union of $B$-orbits. Furthermore, $\overline{\Omega}_{\sigma}\subseteq C_{\sigma}$ for all $\sigma\in\Iu(C_n)$.

It is well-known that $C_w$ is a subvariety of $T_w$ of dimension $\dim C_w=l(w)$ \cite[Chapter 2, Section~2.6]{BileyLakshmibai1}. Let $w\in\Bu(\Phi)$ for $\Phi=B_n$ or $D_n$. Since $\Omega_{w}$ is irreducible, $\overline{\Omega}_{w}$ is. Theorem \ref{theo:dim_Omega} implies $\dim\overline{\Omega}_{w}=\dim\Omega_{w}=l(\sigma)$, so $\overline{\Omega}_{w}$ is an irreducible component of $C_{w}$ of maximal dimension.

\medskip\hypo{Let $\Phi=B_n$ or $D_n$\textup, and $w\in\Bu(\Phi)$. Then the closure of the $B$-orbit $\Omega_{w}$ coincides with the tangent cone $C_w$ to the Schubert variety $\Xu_w$ at the point $p=B/B$.}

\medskip Note that this conjecture implies that if $\sigma\leq_B\tau$, then $\Omega_{\sigma}\subseteq\overline{\Omega}_{\tau}$ for $\sigma,~\tau\in\Bu(\Phi)$.

\end{document}